\documentclass[11pt,final]{article}
\usepackage{amsfonts}
\usepackage{amsmath}
\usepackage{amssymb}
\usepackage{amsthm}
\usepackage{cite}
\usepackage[margin=1.0in]{geometry}
\usepackage{graphicx}
\usepackage{fonts}
\usepackage{color}
\usepackage{algorithm}
\usepackage{algorithmic}
\usepackage{subfigure}

\newcommand{\red}[1]{{\color{red} #1}} 
\renewcommand{\red}[1]{#1} 

\def\bE{{\mathbb{E}}}

\def\dt{{\Delta t}}

\title{On Cellular Automata Models of Traffic Flow with Look-Ahead Potential}

\date{}

\author{
Cory Hauck
\thanks{Computer Science and Mathematics Division,
				Oak Ridge National Laboratory,
				Oak Ridge, TN 37831 USA, (\texttt{hauckc@ornl.gov}).
This author's research was sponsored by the Office of Advanced Scientific
Computing Research and performed at the Oak Ridge National Laboratory,
which is managed by UT-Battelle, LLC under Contract No. DE-AC05-00OR22725. }
\and
Yi Sun
\thanks{Department of Mathematics,
University of South Carolina,
Columbia, SC 29208, (\texttt{yisun@math.sc.edu}).}
\and
Ilya Timofeyev
\thanks{Ilya Timofeyev
Department of Mathematics,
University of Houston,
Houston, TX 77204-3008, (\texttt{ilya@math.uh.edu}).
This author's research was supported in part by the  NSF Grant  
DMS-1109582.
}
}

\begin{document}
\maketitle
\begin{abstract}
We study the statistical properties of a cellular automata model of traffic flow with the look-ahead potential.  
The model defines stochastic rules for the movement of cars on a lattice. 
We analyze the underlying statistical assumptions needed for the derivation
of the coarse-grained model and demonstrate that it is possible to relax some of them to obtain an improved
coarse-grained ODE model. We also demonstrate that spatial correlations play a 
crucial role in the presence of the look-ahead potential and 
propose a simple empirical correction to account for the spatial dependence between
neighboring cells.
\end{abstract}

\tableofcontents

%
%
\section{Introduction}

The derivation of coarse-grained descriptions from microscopic dynamics
has been an active area of research for many decades and a vast amount
of literature exists addressing this issue in various contexts.
In this paper we investigate a particular setup relevant for  
car traffic modeling. Car traffic models can be roughly divided into the following
categories (see review papers \cite{review1,review2,review3} and references therein):
(i) microscopic car-following models that treat cars as particles and
postulate ordinary differential equations (possibly with delay) 
for the car velocity;
(ii) microscopic discrete lattice models where a particular cell 
configuration with values 1 (car is present) and 0 (car is absent) 
combined with explicit rules for movement between cells
are used to represent traffic flow;
(iii) macroscopic partial differential equations models that are typically 
conservation laws relating the car density and flux;
(iv) mesoscopic kinetic models for the velocity distribution, the moments of which gives the macroscopic car density and flux.

While car-following models represent a more realistic setup allowing 
for very detailed interaction rules, lane changing, switching in behavior, etc.,
lattice models are simpler to implement and are more amenable to analytical 
investigation. Therefore, lattice models have been widely used to represent various 
physical phenomena, including car traffic \cite{clq96,nwws98,Sopasakis-Katsoulakis-2006}.
In addition, such cell models are closely related to the Ising-type models and 
cellular automata models and, therefore, a vast literature exists addressing 
various analytical and numerical techniques for models of this type.

One possible use of vehicular traffic models is to utilize filtering techniques to predict the future traffic state.
In this context, particle filtering and Kalman filtering for microscopic models has been developed 
(see e.g. \cite{smh03,wapa05,mibo04}).
Typically, although very detailed rules of interaction can be implemented for microscopic models, 
Monte-Carlo simulations are necessary to estimate statistical properties of
car traffic using these models. Such simulations can be very costly computationally and
using Ensemble Kalman Filtering in conjunction with coarse-grained models is an attractive 
computational alternative \cite{proceed1,proceed2}.

Considerable effort has been devoted to provide the justification of the
macroscopic description including derivations of macroscopic models for
coarse-grained quantities (e.g. car density) from more basic microscopic models.
Recently, a novel look-ahead potential was introduced to model long-range
interactions in prototype lattice model, which were then
coarse-grained to obtain macroscopic descriptions \cite{kmk2003,kms04,kms05}.
In particular, in \cite{Sopasakis-Katsoulakis-2006} the look-ahead potential was
used to model the effect of long-range traffic conditions and a new macroscopic
PDE model with non-local interactions was formally derived. Extensions to multi-lane
traffic have also been developed \cite{duso06,also08}.

In this paper we examine the statistical behavior of the cellular automata (CA)
model introduced in \cite{Sopasakis-Katsoulakis-2006} and the basic underlying
assumptions used in the derivation of the subsequent coarse-grained model. We demonstrate that
while some aspects of the derivation can be improved, assumptions about the approximate independence of the neighboring cells does not hold in general. We also propose a modification of the macroscopic description based on the numerical evidence for the statistical behavior of the
microscopic dynamics. Although the new modified macroscopic model is empirical,
we demonstrate that it can reproduce the results of the stochastic simulations with
remarkable accuracy.

The rest of the paper is organized as follows. In section \ref{sec:meso_old}, we
introduce the CA model, discuss the assumptions about its statistical  behavior, and following
\cite{Sopasakis-Katsoulakis-2006}, outline the derivation of the mesoscopic and
macroscopic models. In section \ref{sec:meso_new} we described an improved
mesoscopic model that is derived by replacing one of the assumptions in \cite{Sopasakis-Katsoulakis-2006} with an exact calculation. In section \ref{sec:num}, we present detailed numerical
experiments to explore the statistical behavior of the microscopic CA model in
various parameter regimes, and in section \ref{sec:hack}, we discuss the
empirical correction to the coarse-grained macroscopic PDE that accounts for the
spatial correlations at the microscopic level. Conclusions are are presented
in section \ref{sec:conclusions}.

%
%
\section{Cellular Automata Model} \label{sec:meso_old}

In this section, we summarize the cellular automata model given in
\cite{Sopasakis-Katsoulakis-2006}.  The model describes a single class of cars
that move in one direction along a single-lane highway.  Possible
multi-lane and multi-class extensions are considered in
\cite{duso06,also08}. The modeled is defined over a one-dimensional
lattice $\cL = \{1,\ldots, N\}$ of $N>1$ evenly spaced cells, and the state of
the system is given by a function $\sigma: \cL \times \bbR^{+} \to \{0,1\}^{N}
$. For $t \in \bbR^{+}$ and $k \in \cL$ , 
\begin{equation}
\sigma_k(t) = \begin{cases}
1\,, & \text{if cell $k$ is occupied at time $t$} \:;\\
0\,, & \text{if cell $k$ is not occupied at time $t$} \:.\\
\end{cases}
\end{equation}
Cars are assumed to move from left to right, and only one car is allowed to occupy a cell at a time. Periodic boundary conditions are imposed so that $\sigma_{mN+k} = \sigma_k$ for any $k \in \cL$ and any integer $m$. 

Transitions in the state of $\sigma$ are the mechanism for modeling car movement.  They obey the rules of an exclusion process \cite{Liggett-1985}:  two lattice sites exchange values in each transition and cars may not move into occupied cells.  In addition cars are only allowed to move one cell to the right.  Thus the only possible configuration changes are of the form
\begin{equation}
\label{eq:config}
\left\{ \sigma_k(t)=1,\sigma_{k+1}(t)=0 \right\}
\to
\left\{ \sigma_k(t+\Delta t)=0,\sigma_{k+1}(t + \Delta t) = 1 \right\} \:.
\end{equation}
In the simplest case, the transition rate $c_0>0$ is a constant.  Thus the configuration change in \eqref{eq:config}
occurs with probability $c_0 \Delta t + o(\dt)$ for small
$\Delta t$, and in the absence of other forces, cars move with velocity
$v_0 = h c_0$ on average, where $h$ is the cell width.  As the car velocity
does not depend on the cell width, $c_0$ must scale with $h^{-1}$.

In \cite{Sopasakis-Katsoulakis-2006}, the transition rate is given by an Arrhenius-type
formula \cite{Laidler-1987} with a look-ahead potential $J_k$ that depends on values of
$\sigma_l$ for $l > k$. In this case, the transition rates becomes $c_0
e^{-\beta J_k}$, where
\begin{equation}
\label{Jk}
J_k(t) = \frac{1}{M}\sum_{i=1}^{M} \sigma_{k+i+1}(t),
\end{equation}
$\beta > 0$ is a parameter describing the strength of the look-ahead
interactions, and $M$ is the number of cars to the right of cell $k$ which
affect velocity of the car in cell $k$. The term $e^{-\beta J_k}$ plays the role
of a slowdown factor when the forward car density is high, i.e, when the road is
congested. It is also possible to introduce weights into \eqref{Jk} so that cars
nearby have a more pronounced effect than those that are further away.

The microscopic model described above can be formulated as a continuous-time
Markov chain and the generator of this process can be computed explicitly.
It is defined as 
\begin{equation}
(A \psi)(t) = \lim_{\Delta t \to 0} \frac{\bE \psi(\sigma(t+\Delta t)) - \psi(\sigma(t))}{\Delta t} \:,
\end{equation}
where $\psi$ is any test function and the expectation is 
taken over all possible transitions between time $t$ and $t+\Delta t$.
For an arbitrary test function $\psi$, the generator is given by
\begin{equation}
\label{eq:generator}
A \psi = \sum_{k \in \cL}
c_0 e^{-\beta J_k} \sigma_k [1-\sigma_{k+1}]
\left[ \psi(\sigma^{k,k+1}) - \psi(\sigma) \right] \:,
\end{equation}
where $\sigma^{k,k+1}$ denotes the lattice configuration $\sigma$ after an exchange between 
the cells $k$ and $k+1$:
\begin{equation}\
\sigma^{k,k+1}_l = 
\begin{cases} 
\sigma_l \:, & l \ne k,k+1 \:,\\
\sigma_{k+1}\:, & l = k \:,\\
\sigma_{k}\:, & l = k+1\:.
\end{cases}
\end{equation}
The terms $\sigma_k$ and  $[1-\sigma_{k+1}]$ in \eqref{eq:generator}
ensure that only configuration changes of the form given in \eqref{eq:config}
contribute to the sum, i.e., that before a transition can occur, there must be a car in cell $k$ and cell $k+1$ must be empty.

A coarse grain model is derived by first computing the action of the generator on the test function 
$\psi(\sigma) = \sigma_l$. In this case, the only indices that contribute to the sum in \eqref{eq:generator} are $k=l-1$ and $k=l$.
Thus, $A \sigma_l$ becomes
\begin{align}
A \sigma_l &= c_0 e^{-\beta J_{l-1}} \sigma_{l-1} [1-\sigma_{l}]
\left[\sigma^{l-1,l}_l - \sigma_l \right] + 
c_0 e^{-\beta J_l} \sigma_l [1-\sigma_{l+1}]
\left[\sigma^{l,l+1}_l - \sigma_l \right] \nonumber \\
&=  c_0 e^{-\beta J_{l-1}} \sigma_{l-1} [1-\sigma_{l}]
- c_0 e^{-\beta J_l} \sigma_l [1-\sigma_{l+1}].
\end{align}
From the definition of the
generator,
\begin{equation}
\label{evol}
\frac{d}{dt} \bbE (\psi) = \bbE (A \psi) \:,
\end{equation}
where $\bbE$ is the expectation operator with respect to the probability
measure associated with the stochastic process.
Thus the evolution equation for $\rho_k := \bbE (\sigma_k$) is given by
\begin{equation}
\label{den1}
\frac{d}{dt} \rho_k =  \bbE \left(
c_0 e^{-\beta J_{k-1}} \sigma_{k-1} (1-\sigma_{k}) -
c_0 e^{-\beta J_k} \sigma_k (1-\sigma_{k+1})
\right).
\end{equation}

A mesoscopic model for $\rho = \{\rho_1,\ldots,\rho^N\}$ requires a closure
that approximates the right-hand side of \eqref{den1} by a function of $\rho$.
In
\cite{Sopasakis-Katsoulakis-2006}, a closure was found under two assumptions:
\begin{itemize}
\item[\bf A1.] The look-ahead interactions are weak so that
$\bbE
e^{-\beta J_k} \approx e^{-\beta \bbE J_k}$.
\item[\bf A2.] The probability measure on $\sigma$ is approximately a
product measure, i.e., $f(\sigma_k = \delta_k, \sigma_l = \delta_l) =
f(\sigma_k = \delta_k) f(\sigma_l = \delta_l)$, where $f$ is the probability
density.
\end{itemize}
Note that Assumption A2 implies, among other things, that $\sigma_k$ and
$\sigma_l$ are uncorrelated, i.e.
\begin{equation}
   \bbE\left(\sigma_k(t) \sigma_{l}(t)\right)
\approx
\bbE(\sigma_k(t)) \, \bbE\left(\sigma_{l}(t)\right) \equiv \rho_k \rho_l  
\end{equation} for $k \ne l$.  While this weaker condition is sufficient for
closure without the look-ahead, independence of higher-order moments is needed
for closure with the look-ahead.

Based on the assumptions above, the resulting equation for the density is
\begin{equation}
\label{meso1}
\frac{d}{dt} \rho_k =  
c_0 e^{-\beta I_{k-1}} \rho_{k-1} (1-\rho_{k}) -
c_0 e^{-\beta I_k} \rho_k (1-\rho_{k+1}) \:,
\end{equation}
where 
$$
I_k = \frac{1}{M} \sum_{i=1}^{M} \rho_{k+i+1}
$$
and $\rho_{mN+j} = \rho_l$ for all integers $j$.

The validity of the two assumptions above depends on the strength of the
look-ahead potential. Assumption A1 implies that $\beta$ is small, in which case
the influence of the look-ahead potential $J_k$ is weak. Assumption A2 means
that $\sigma_k(t)$ and $\sigma_{k+i}(t)$ are independent for all $i$, including
the nearest neighbors $i=1,2$.
\red{
This is somewhat unrealistic even without the look-ahead potential, since the exclusion principle
couples the neighboring cells. Nevertheless, our simulations show that without the look-ahead potential
the coupling is weak in most situations and the approximate independence assumption is quite plausible.
On the other hand, in the presence of a strong look-ahead potential the neighboring cells 
become very highly anti-correlated. This can be understood as follows:  If the
length of the look-ahead potential is $M$, then a car in position $k$ is
very likely to be trailed by $M$ zeros
since the probability to move for the car in the position $k-M-1$ is $c_0 e^{-\beta/M} \Delta t \ll c_0 \Delta t$ for large $\beta$. Therefore, the probability to move is extremely small, and a car in the position $k-M-1$ is very likely to ``wait'' for the car in the position $k$ to move.  Positive correlations may also occur in some situations with and without the look-ahead potential.
We discuss such situations in the numerical experiments in section
\ref{sec:num}.} 
In summary,
assumptions A1 and A2 are plausible if the look-ahead interactions are weak
($\beta \ll 1$), but are unlikely to hold if the look-ahead interactions are
stronger ($\beta = O(1)$).

A first-order PDE model can be formally derived from the mesoscopic model in the limit $N
\to \infty$. Let the look ahead distance $L$ and domain length $D$ be fixed;
set $h = D/N$ and $M = [L/D]N$.  For $k \in \{1,\ldots,N\}$, let $x_k =(k-0.5)h$
denote
the center of cell $k$, and let
$\bar{\rho}^N$ be a smooth function of $x$ and $t$ such that
$\bar{\rho}^N(x_k,t)
= \rho_k(t)$.
Define
\begin{equation}
  F^N(x,t) = e^{-\beta \bar{I}^N(x,t)} \bar{\rho}^N(x-h,t) [1-\bar{\rho}^N(x,t)]
\:.
\end{equation}
where $\bar{I}^N$ is a smooth interpolant of $I_k$:
\begin{equation}
  \bar{I}^N(x) = \frac{1}{M} \sum_{i=1}^M \bar{\rho}^N(x+(i+1)h)
  = \frac{1}{L}\int_0^L \bar{\rho}^N(x+y,t) dy + O(h) \:.        
\end{equation}
 Then with a Taylor expansion of
$F^N$ around $x_k$, \eqref{meso1} becomes
\begin{equation}
\frac{\partial}{\partial t} \bar{\rho}^N(x_k,t)
= c_0 ( F^N(x_{k-1},t) - F^N(x_k,t) )
= v_0 \frac{\partial}{\partial x}F^N(x_k,t) + O(h) \:.
\end{equation}
One can then formally pass to the limit $N \to \infty$ while $L$ and $D$
remain fixed.  For each fixed $x \in (0,D)$, let $k = [x/h]$ so that $x_k \to
x$ 
as $h \to 0$ and $\bar{\rho} = \lim_{N \to \infty} \bar{\rho}^N$ satisfies
\begin{equation}
\frac{\partial}{\partial t} \bar{\rho}(x,t)
+ \frac{\partial}{\partial x} \phi(\bar{\rho}) = 0 \:,
\label{eq:cons_law}
\end{equation}
where the non-local flux $\phi$ is given by
\begin{equation}
\phi(\bar{\rho})(x,t) = v_0 e^{-\frac{\beta}{L} \int_0^L \bar{\rho}(x+y,t) dy}
\bar{\rho}(x,t)
(1-\bar{\rho}(x,t)) \:.
\end{equation}

%
%
\section{A New Mesoscopic Model} \label{sec:meso_new}

In this section we show the expectation of $e^{-\beta J_k}$ can be computed
exactly so that assumption A1 can be removed. This leads to a slightly improved
mesoscopic model for the car density. (However, as in the previous mesoscopic
model, assumption A2 is still required to close products.)  In particular, A2
allows us to approximate
\begin{equation}
\bE \left( e^{-\beta J_k} \right)
= \bE \left( \prod_{i=1}^{M} e^{-\beta' \sigma_{k+i+1}} \right)
\simeq \prod_{i=1}^{M} \bbE \left( e^{-\beta' \sigma_{k+i+1}} \right) \:.
\label{prod_meas}
\end{equation}
where $\beta' =  {\beta}/{M}$.
We then use the fact that for $1 \leq l \leq N$ and for any positive integer
$m$, $(\sigma_l)^m = \sigma_l$.  Simple algebra gives
\begin{eqnarray}
e^{-\beta'\sigma_l}=
\sum_{m=0}^{\infty} \frac{(-\beta')^m (\sigma_l)^{m}}{m!}
= 1 + \sigma_l \sum_{m=1}^{\infty} \frac{(-\beta')^m}{m!}
= 1 + \sigma_l ( e^{-\beta'} - 1) \:.
\end{eqnarray}
Thus the expectation of
$e^{-\beta' \sigma_l}$ is
\begin{equation}
\bbE (e^{-\beta\sigma_l} )
= 1 + \bbE \left(\sigma_l \right)( e^{-\beta'} - 1)
= 1 + \rho_l( e^{-\beta'} - 1)  \:.
\label{expsigma} 
\end{equation}
Substituting  \eqref{expsigma} into the right hand side of \eqref{prod_meas}
gives
\begin{equation}
\label{expjk}
\bE \left( e^{-\beta J_k} \right)
= \prod_{i=1}^{M} \left[ 1 + \rho_{k+i+1} \left( e^{-\beta'} - 1 \right)
\right] \:,
\end{equation}
This expression can be used to derive a new mesoscopic model for the
car density $\rho_k(t)$:

\begin{align}
\frac{d}{dt} \rho_k =&
c_0 \rho_{k-1} (1-\rho_k) 
\prod_{i=1}^{M} \left[ 1 + \rho_{k+i} \left( e^{-\beta'} - 1 \right) \right]
\nonumber \\&-
c_0 \rho_k (1-\rho_{k+1}) 
\prod_{i=1}^{M} \left[ 1 + \rho_{k+i+1} \left( e^{-\beta'} - 1 \right)
\right]\:,
\label{meso2}
\end{align}
which can be rewritten as
\begin{align}
\frac{d}{dt} \rho_k =&
c_0 \rho_{k-1} (1-\rho_k) 
e^{\sum_{i=1}^{M} \log \left[ 1 + \rho_{k+i} \left( e^{-\beta'} - 1 \right)
\right]}
\nonumber \\
 &-c_0 \rho_k (1-\rho_{k+1})
e^{\sum_{i=1}^{M} \log \left[ 1 + \rho_{k+i+1} \left( e^{-\beta'} - 1 \right)
\right]}
\label{meso22} \:.
\end{align}

One can easily recover the model in \eqref{meso1} by expanding
the exponential and logarithm in \eqref{meso22} to leader order in $\beta'$. 
Indeed, simple Taylor expansions give
\begin{equation}
\log \left[ 1 + \rho_{k+i+1} \left( e^{-\beta'} - 1 \right)
\right]= -\beta'\rho_{k+i+1} + O(\beta'^2)
\end{equation}
In particular, for fixed $\beta$ and large $M$,
\begin{equation}
\sum_{i=1}^M
\log \left[ 1 + \rho_{k+i+1} \left( e^{-\beta'} - 1 \right) \right]
= -\beta I_k + O(M^{-1}) 
\end{equation}
and since $M^{-1}$ = O($h$) as $N \to \infty$, the new mesoscopic model has the
same PDE limit as the old one.

%
%
\section{Numerical Experiments}
\label{sec:num}
In this section we perform a detailed comparison between the stochastic traffic
model and the two mesoscopic models \eqref{meso1} and \eqref{meso22} in various
parameter regimes. The new mesoscopic model in \eqref{meso22} provides a better
front tracking in some regimes. Nevertheless, as discussed earlier,
it still suffers from the same limitation as the old model, namely that
look-ahead interactions should be weak so that the measure associate with
$\sigma$ is nearly a product measure. When the
look-ahead interactions are
weak ($\beta$ is small), both mesoscopic models perform similarly.

\subsection{The Stochastic Algorithm} We use the Metropolis algorithm
\cite{MetropolisAlgorithm} to simulate the stochastic cellular automata model.%
\footnote{
The Metropolis algorithm was cross-validated with a Kinetic
Monte-Carlo algorithm \cite{voter-2007}.  Several
comparisons were made to ensure statistically accurate results.}
Given a
time-step $\delta
t$, we advance the configuration of the lattice during each time-step according
to the following algorithm. If a car is present at the position $k$ (i.e.
$\sigma_k(t)=1$) and there is no car at the position $k+1$ (i.e.
$\sigma_{k+1}(t)=0$), then the car at the position $k$ can advance with
probability $c_0 e^{-\beta J_k}\delta t$, where $J_k$ is the
look-ahead potential in \eqref{Jk}. Simple pseudo-code for a single time step is
given below.

\begin{algorithm}
\caption{Metropolis Algorithm (one time step)}
\label{alg:metropolis algorithm}
\begin{algorithmic}
\REQUIRE $N$, $M$, $t$, $\beta$,  $\delta t$, and $\sigma_k(t)$ for
$k=1:N$
\FOR{ $k=1:N$}
  \IF{$\sigma_k(t) == 1$ \textbf{AND} $\sigma_{k+1}(t) == 0$}
    \STATE  $J_k = \frac{1}{M}\sum_{i=1}^{M} \sigma_{k+i+1}(t)$
    \STATE  $p =$ unif(0,1) \\
    \IF{$p < C_0 e^{-\beta J_k} \delta t$}
      \STATE  $\sigma_k(t+\delta t) = 0$
      \STATE  $\sigma_{k+1}(t+\delta t) = 1$
    \ENDIF
  \ENDIF
\ENDFOR
\end{algorithmic}
\end{algorithm}

The configuration array $\sigma_k(t+\delta t)$ is stored separately from
$\sigma_k(t)$. Thus, changes in the configuration in one part of the lattice do
not affect the calculation of $J_k$ for other indexes in the lattice.

All numerical simulations presented in this paper are ensemble Monte-Carlo
simulations and expectations plotted in all Figures are ensemble averages.
Given the initial conditions \eqref{ic1}, we generate
$n=5000$ different realizations of $\sigma$ and estimate the density by the
formula
\begin{equation}
\rho(k,t) \approx \frac{1}{n} \sum\limits_{p=1}^{n} \sigma_k^{(p)}(t) \:,
\end{equation}
where the integer $p$ is the realization index. During the course of the simulations, various statistics are recorded in order
to make the comparisons presented below.

We perform all experiments using the deterministic ``red light'' initial
condition
\begin{equation}
\sigma_k(0) =
\begin{cases}
 1 & 20 \le k \le K, \\
 0 & K < k \le N,
\end{cases}
\label{ic1}
\end{equation}
with $K=60$.  This condition allows us to examine the behavior of the leading
rarefaction wave.  Many other configurations are possible including those which
lead to shock waves. Following \cite{Sopasakis-Katsoulakis-2006}, we set the
cell size to 22 feet. Other parameters in the simulation are $c_0 = 1/0.23 =
4.3478$ (corresponds to approximately $65$ miles/hour) and $N=700$.  The
boundary conditions are periodic, but $N$ is set large enough to ensure that
they do not affect the simulations.

\subsection{Assumption A1:  Expectation of the Exponential}

\begin{center}
\begin{figure}[H]
\subfigure[$\beta=1$, $M=1$]{
\includegraphics[width = 0.45\textwidth]{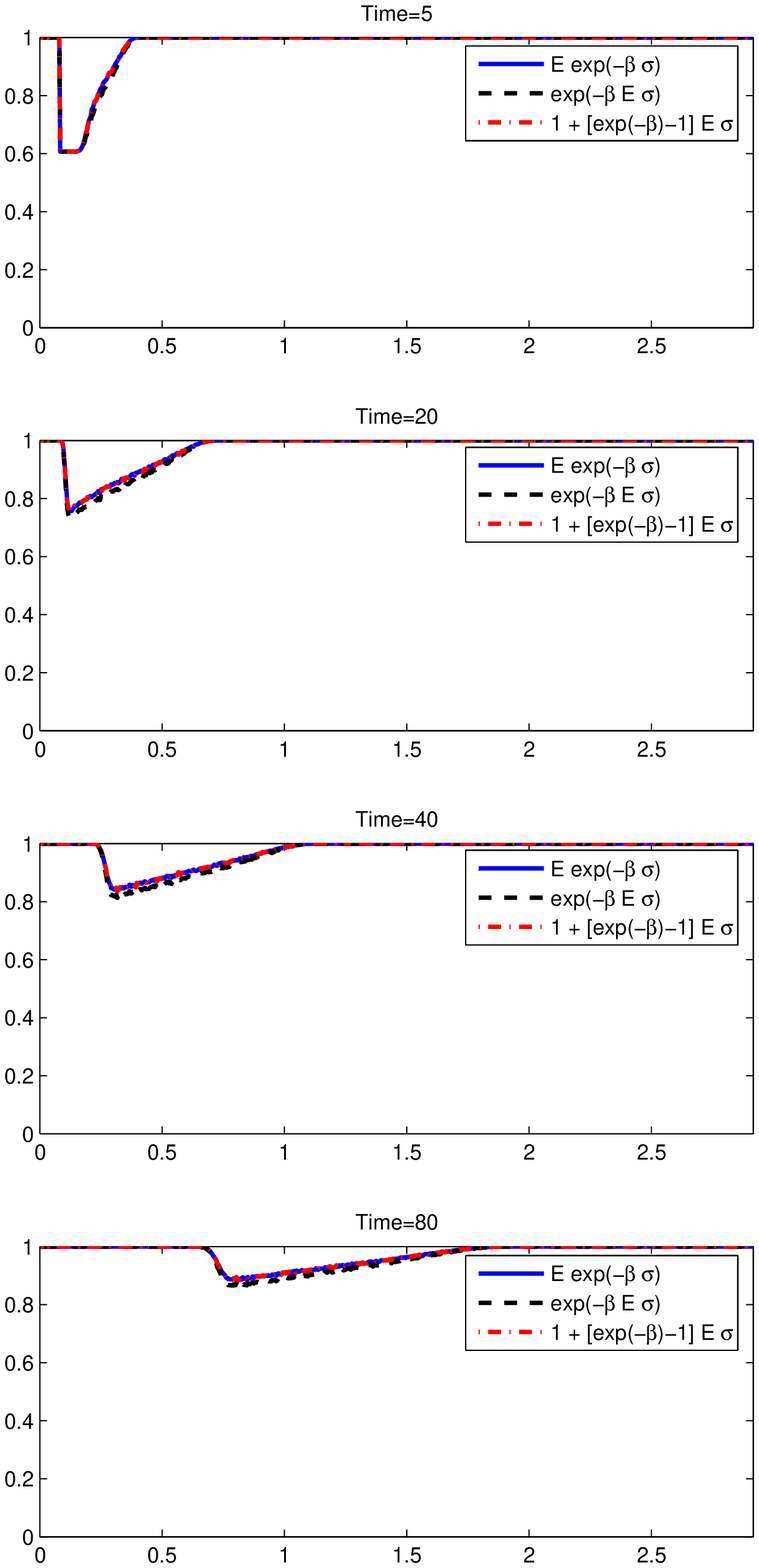}}
\subfigure[$\beta=3$, $M=5$]{
\includegraphics[width = 0.45\textwidth]{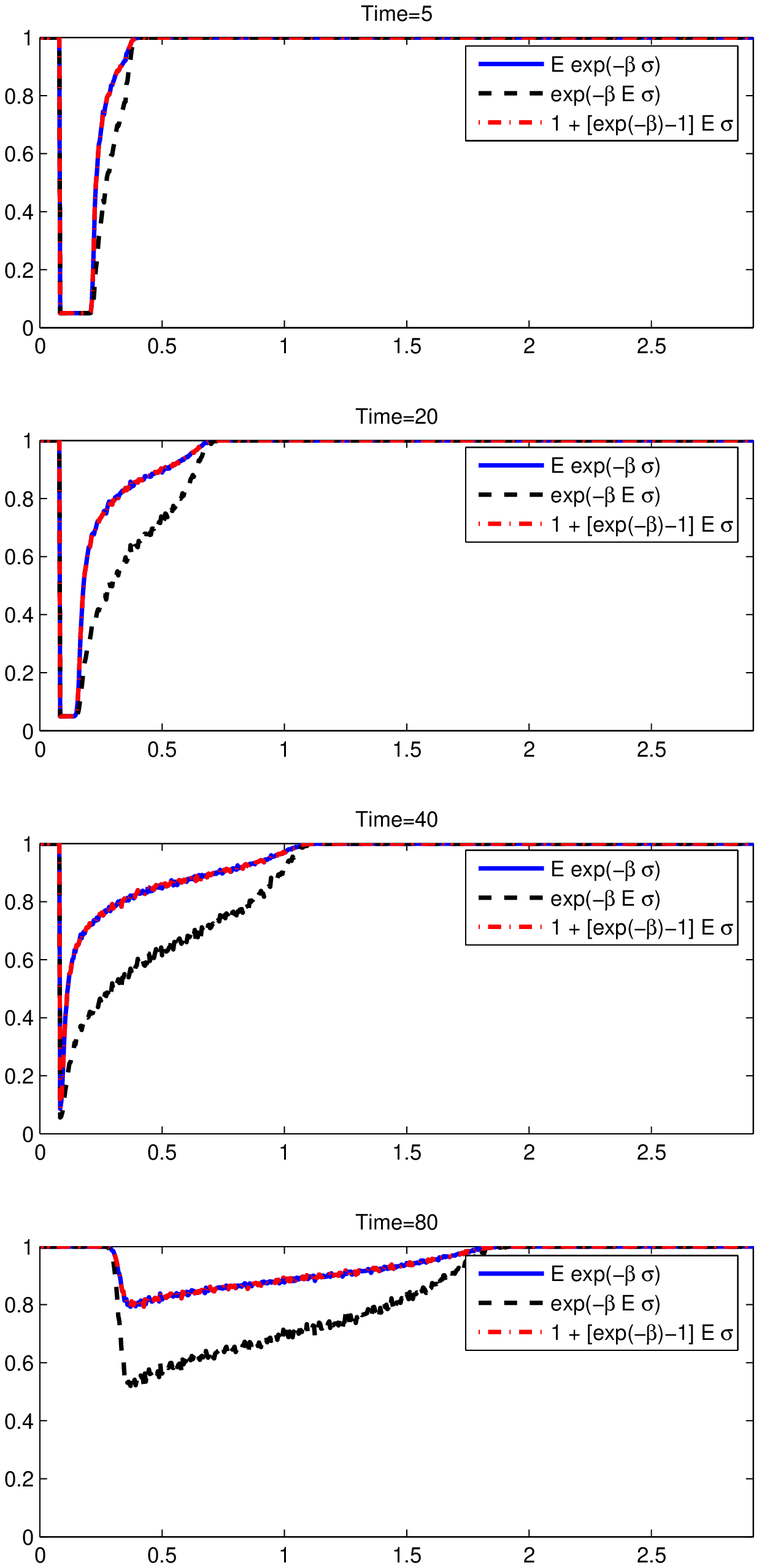}}

\caption{Comparison of $\bE(e^{-\beta \sigma_k})$ with the closure
$\bE(e^{-\beta \sigma_k}) \approx e^{-\beta \bE \sigma_k}$ 
the closure in \eqref{expsigma}. Note that the solid blue and the red line overlap completely.}
\label{figexpsig}
\end{figure}
\end{center}

We first check that the new mesoscopic model removes the error in the
closure assumption A1, i.e. that $\bE(e^{-\beta \sigma_k}) \approx e^{-\beta\bE
\sigma_k}$.  Two representative cases are given in
Figures \ref{figexpsig}: one with $\beta=0.5$ and $M=1$ and a second with
$\beta = 3$ and $M=5$.  For reference, we include the formula in
\eqref{expsigma}, which is used in the
new model. As expected, this formula is exact for all values of $\beta$ and
$M$.  The dependence of the results on $\beta$ is generic:  Assumption A1 is
nearly valid for values of $\beta < 1$, since in this case the powers of $\beta$
quickly decay in the expansion of the
exponential.  However, the accuracy of the approximation quickly decreases for
larger values of $\beta$.  Many other parameter values have been tested to
confirm this.  Roughly speaking, the validity of Assumption 1 appears to be largely independent of
the size of $M$.

\subsection{Front Tracking}

\begin{center}
\begin{figure}[H]
\subfigure[$\beta=0.5, M=5$]{
\includegraphics[width=0.45\textwidth]{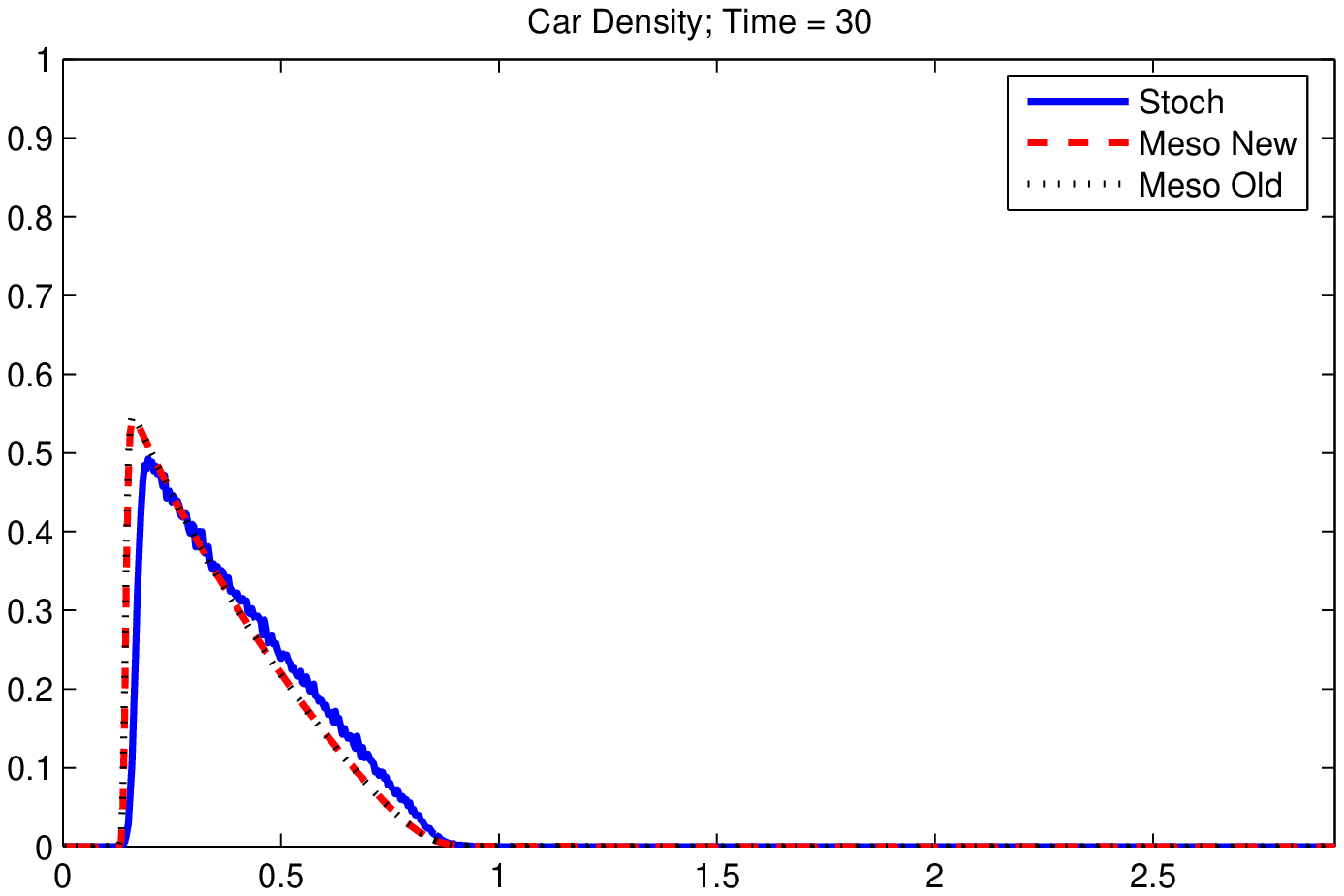}}
\subfigure[$\beta=0.5, M=5$]{
\includegraphics[width=0.45\textwidth]{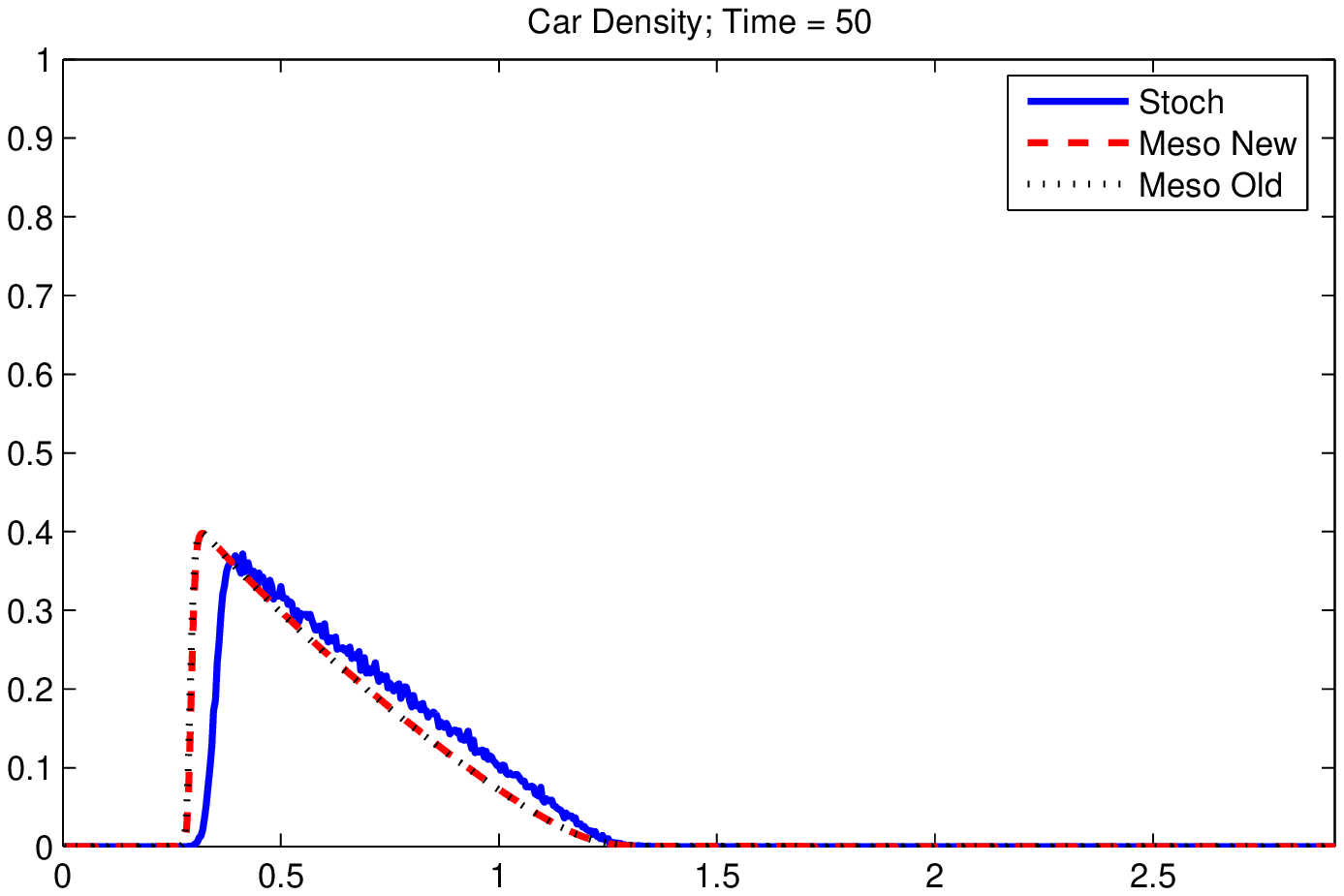}}
\\~\\
\subfigure[$\beta=3, M=5$]{
\includegraphics[width=0.45\textwidth]{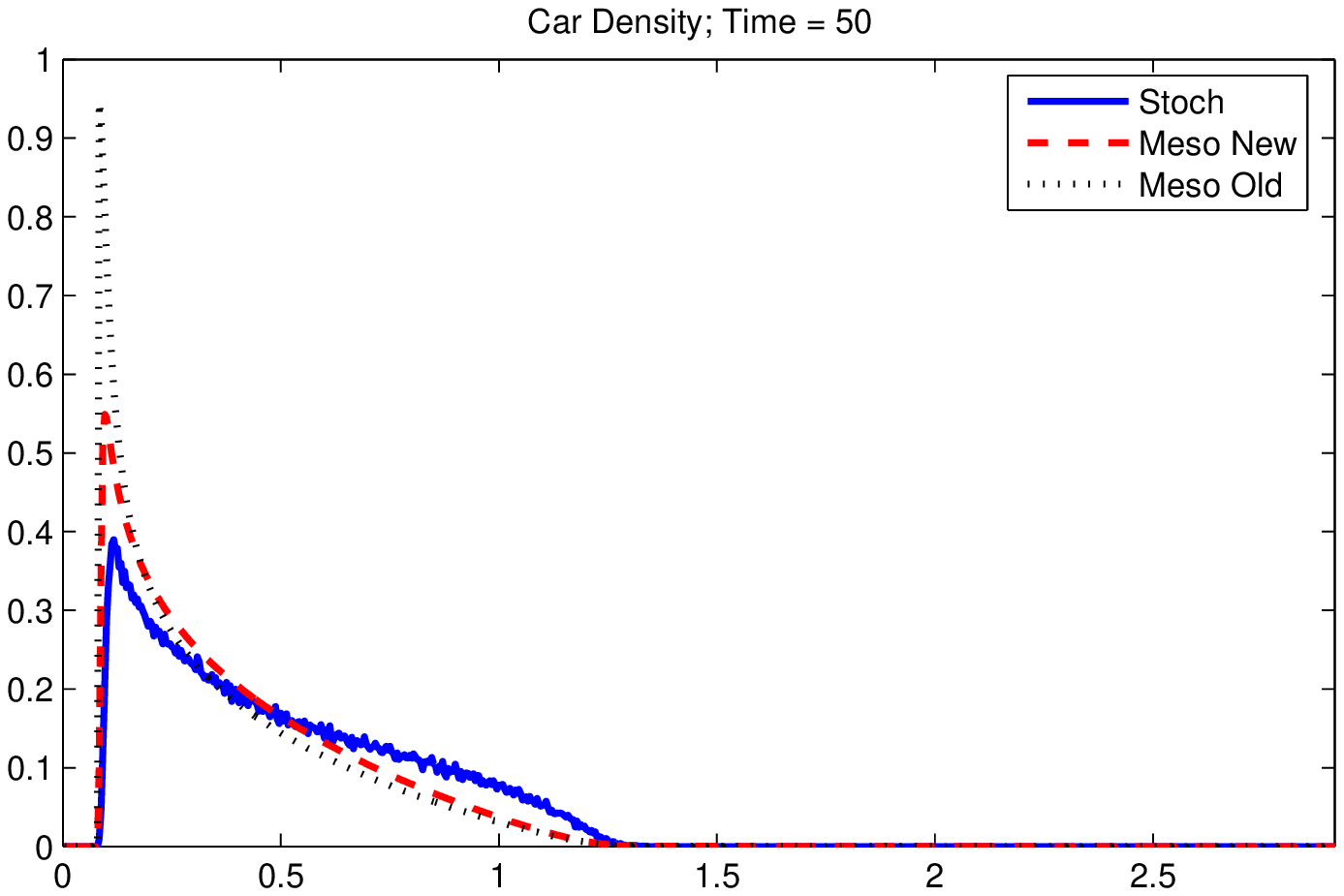}}
\subfigure[$\beta=3, M=5$]{
\includegraphics[width=0.45\textwidth]{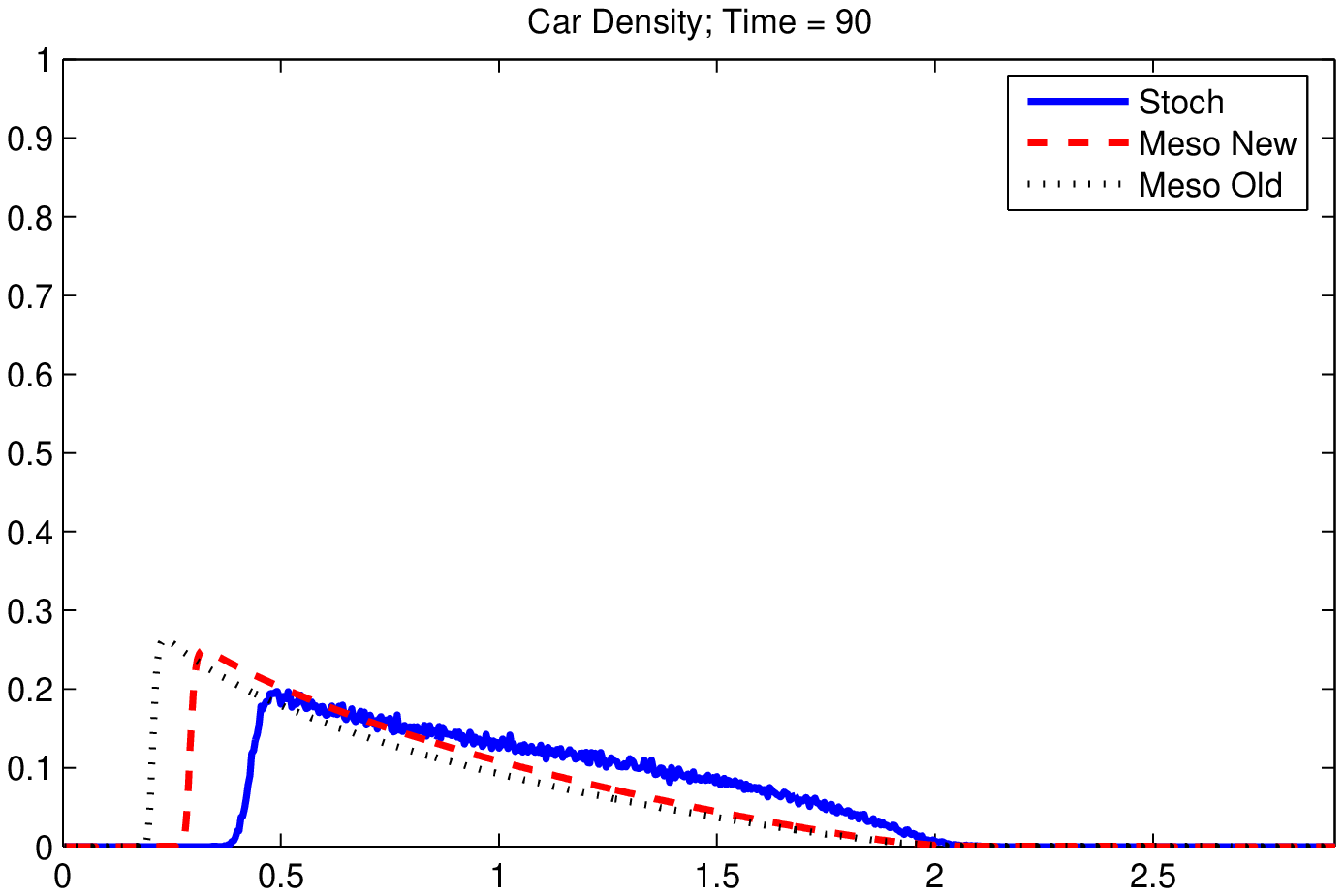}}
\\~\\
\subfigure[$\beta=5, M=5$]{
\includegraphics[width=0.45\textwidth]{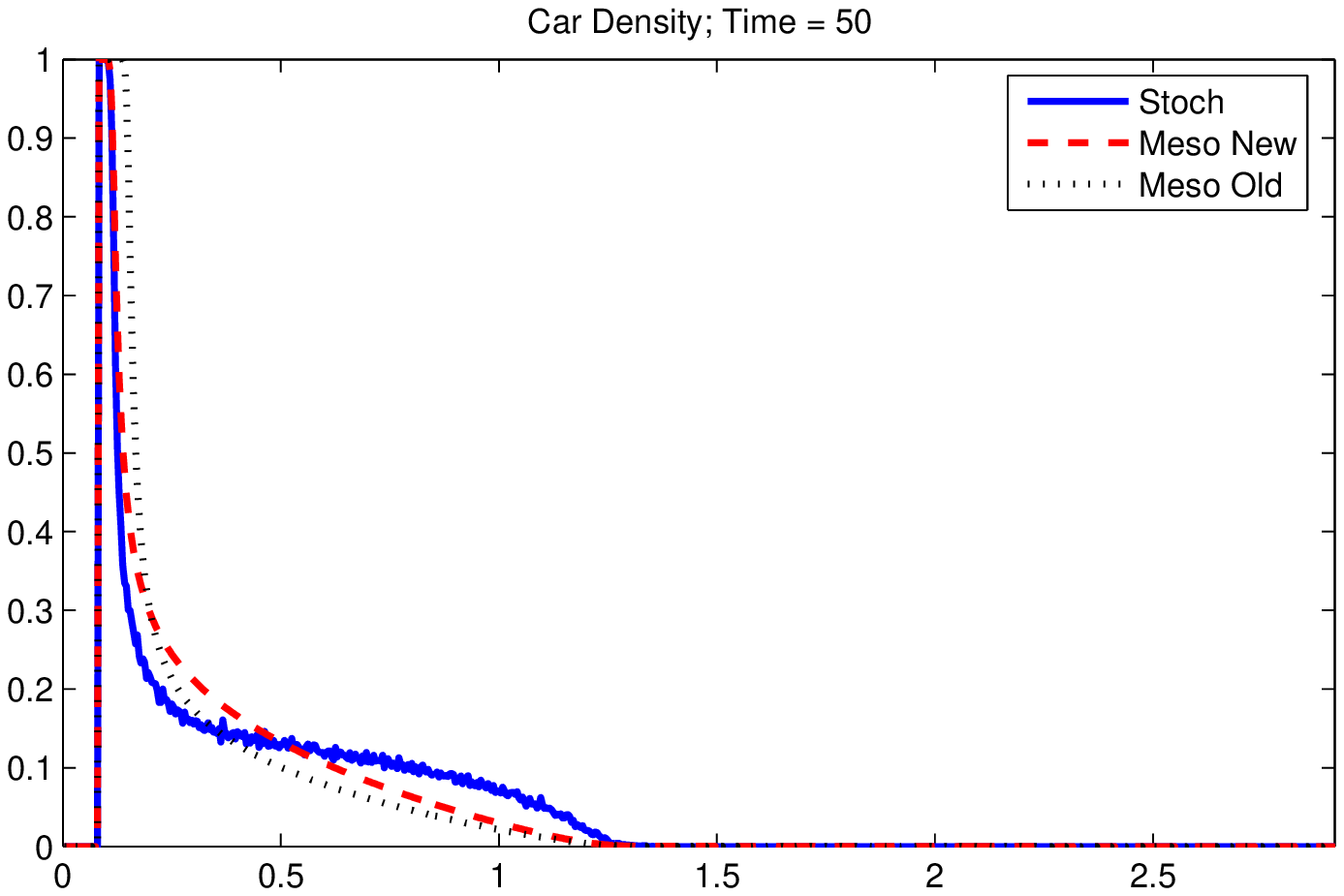}}
\subfigure[$\beta=5, M=5$]{
\includegraphics[width=0.45\textwidth]{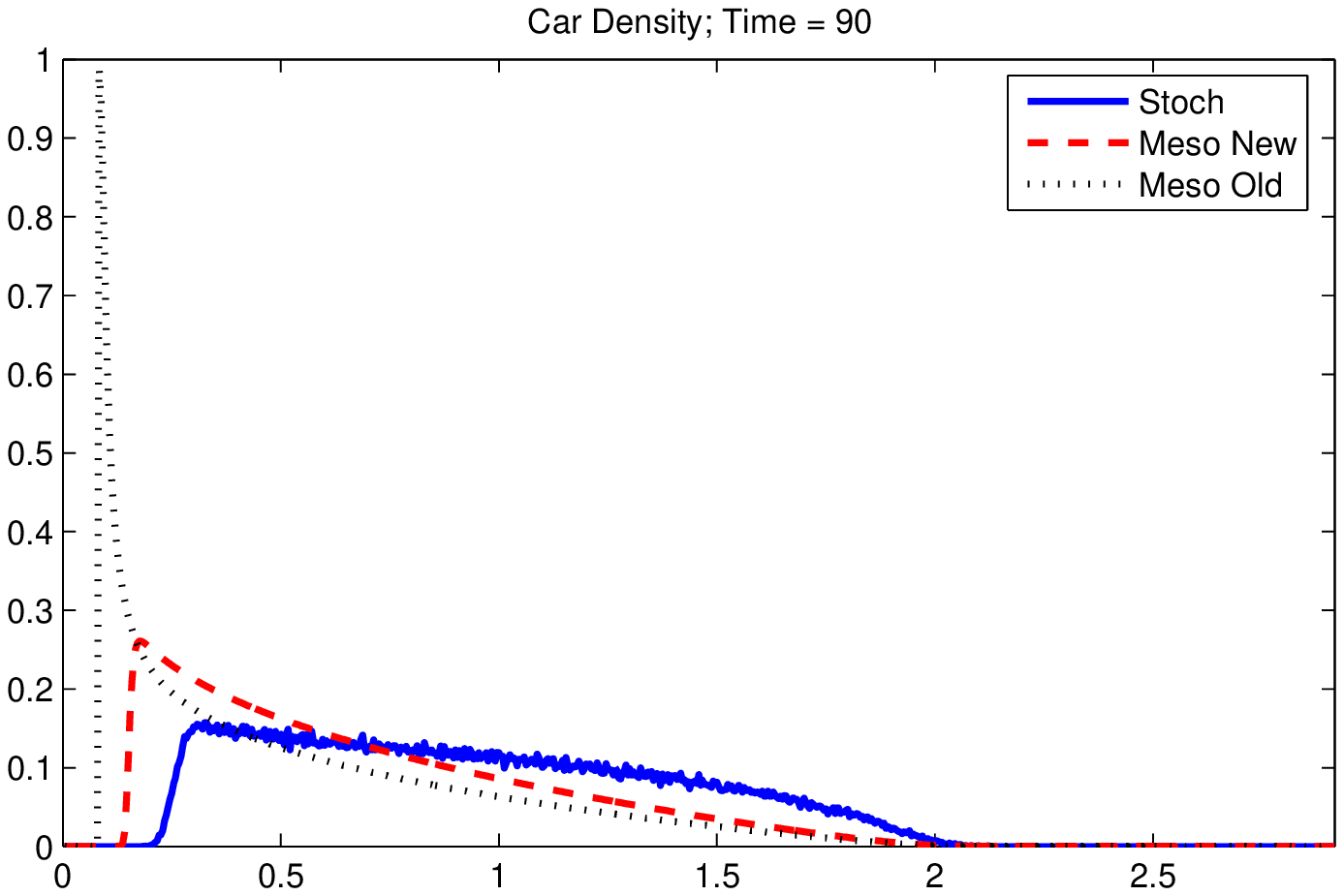}}
\caption{Comparison of Simulations of the Stochastic and Mesoscopic Models with
the initial condition \eqref{ic1}. Note that for $\beta = 0.5$, the numerical results from the
mesoscopic models
\eqref{meso1} and \eqref{meso22} overlap in this regime.}
\label{fig1}
\end{figure}
\end{center}

In our next set of experiments, we investigate how well each of the mesoscopic
models tracks the rarefaction front.  We consider the look-ahead distance $M=5$
and vary the strength of the look-ahead interaction $\beta$. A comparison of the
stochastic simulation and both mesoscopic models is given in Figure \ref{fig1}.
In all cases, the bulk of the traffic proceeds faster in the stochastic
simulation and the trailing front is less steep, with the new model always
outperforming the old one.  As expected, the difference between the stochastic
and mesoscopic models increases as $\beta$ is increased. Indeed, for larger
values of $\beta$, there is a considerable discrepancy between the stochastic
simulations and both mesoscopic models.  However, front tracking with the new
mesoscopic model \eqref{meso22} is still considerably better than with the old
model \eqref{meso1}.

\subsection{Assumption A2:  Correlations}

To assess the source of the discrepancy between the microscopic model and the
new mesoscopic model in \eqref{meso22}, we
investigate numerically the validity of assumption A2---that the measure
associated with $\sigma$ is (nearly) a product measure.  To do this, we examine
for small $i \ge 1$
the size of the sample Pearson correlation coefficient \cite{Devore-2011}:%
\footnote{Note that a zero correlation coefficient is a necessary but not
sufficient condition to have a product measure.
}
\begin{align}
r_{k,k+i}(t) &= \frac{1}{n-1} \sum ^n _{p=1} \left( \frac{\sigma^{(p)}_{k}(t) -
\bar{\sigma}_{k}(t)}{s_k(t)} \right) \left( \frac{\sigma^{(p)}_{k+i}(t) -
\bar{\sigma}^{(p)}_{k+i}(t)}{s^{(p)}_{k+i}(t)} \right) \nonumber \\
&= \frac{\sum ^n _{p=1}(\sigma^{(p)}_{k}(t) -
\bar{\sigma}_{k}(t))(\sigma^{(p)}_{k+i}(t) - \bar{\sigma}^{(p)}_{k+i}(t))}{\sqrt{\sum ^n
_{p=1}(\sigma^{(p)}_{k}(t) - \bar{\sigma}_{k}(t))^2} \sqrt{\sum ^n
_{p=1}(\sigma^{(p)}_{k+i}(t) - \bar{\sigma}_{k+i}(t))^2}} \:,
\end{align}
where $s_k(t)$ is the standard deviation of the stochastic process at cell $k$
and time $t$, bars denote sample means, and $p$ indexes the samples.

\begin{center}
\begin{figure}[H]
\centerline{
\includegraphics[width=1.5in,height = 6.0in]{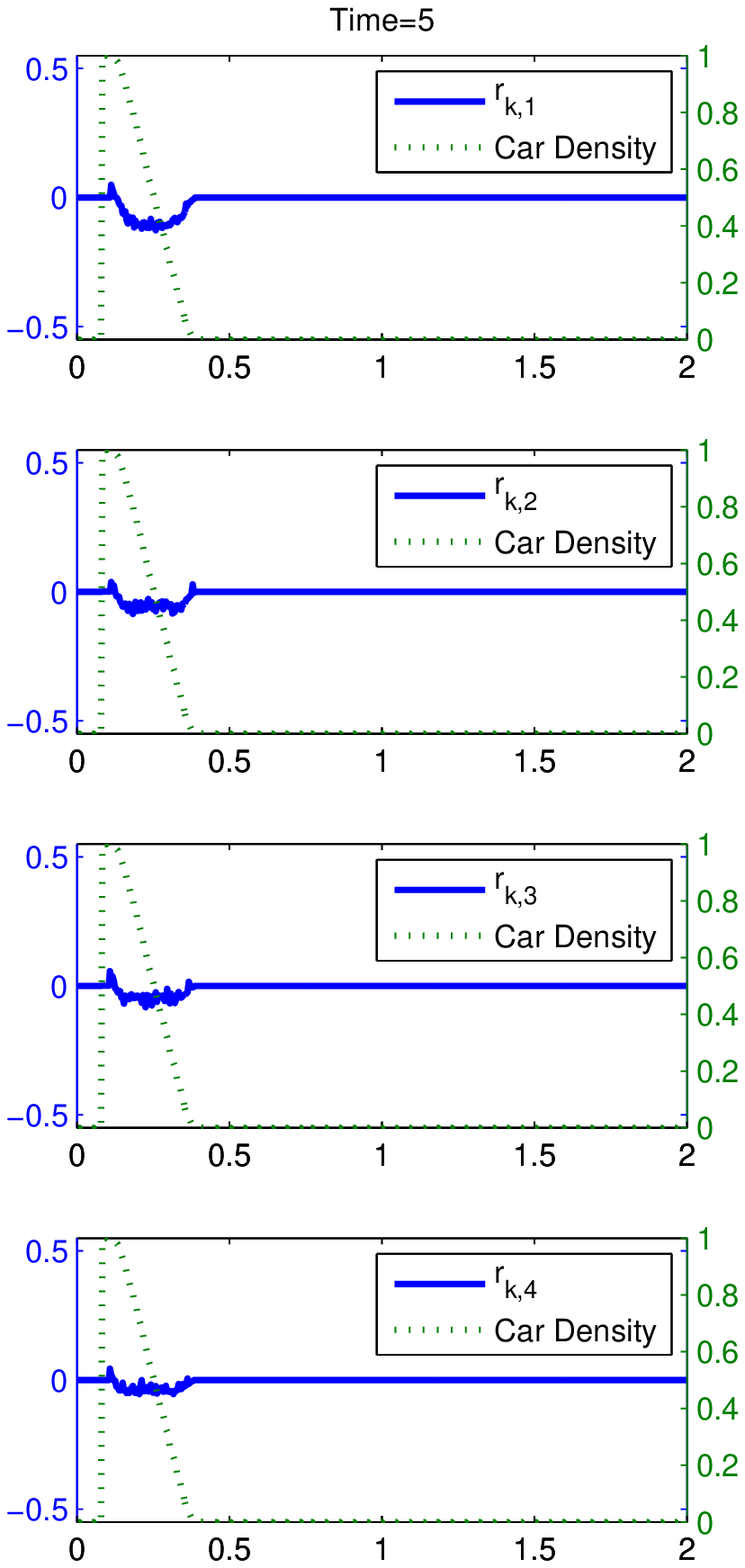}
\includegraphics[width=1.5in,height = 6.0in]{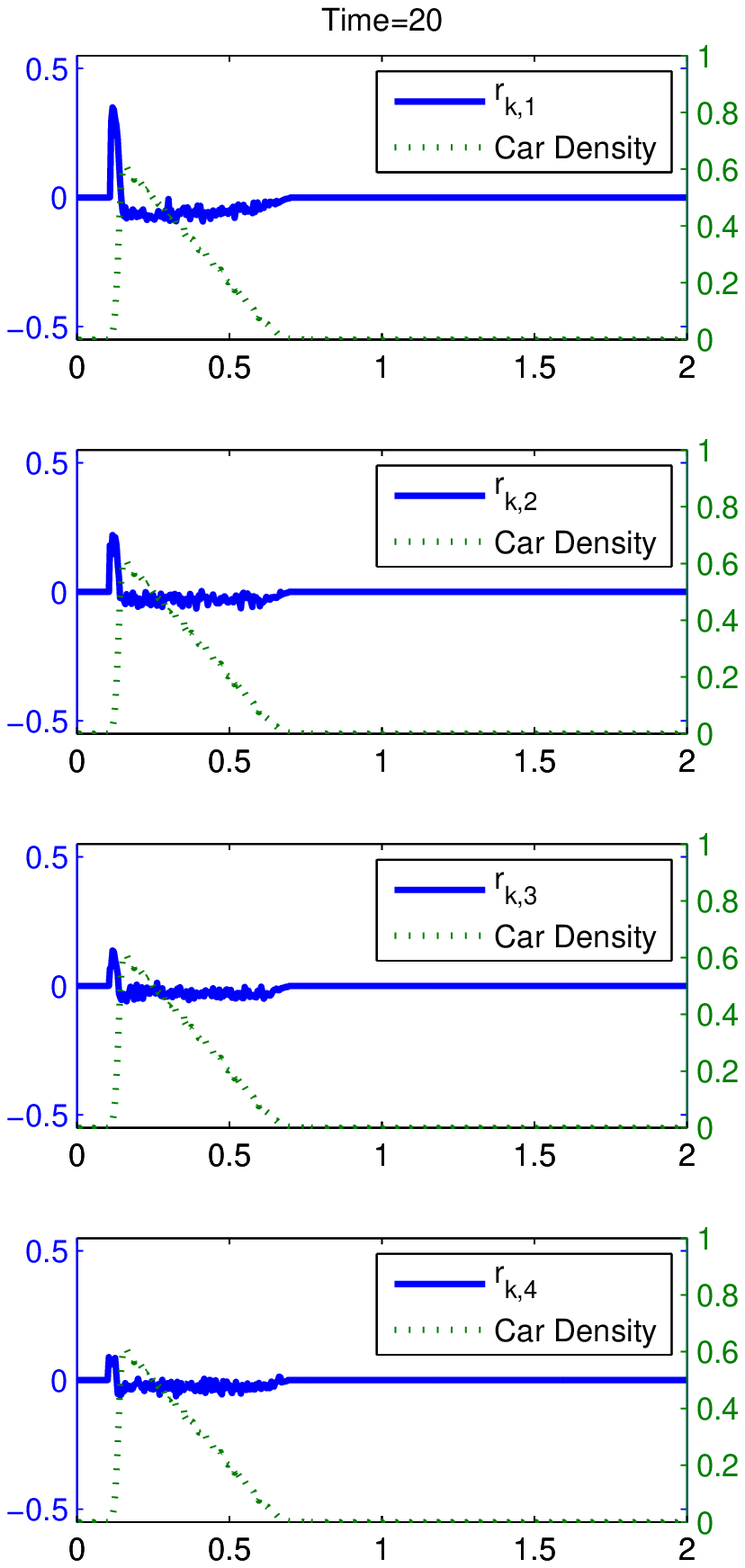}
\includegraphics[width=1.5in,height = 6.0in]{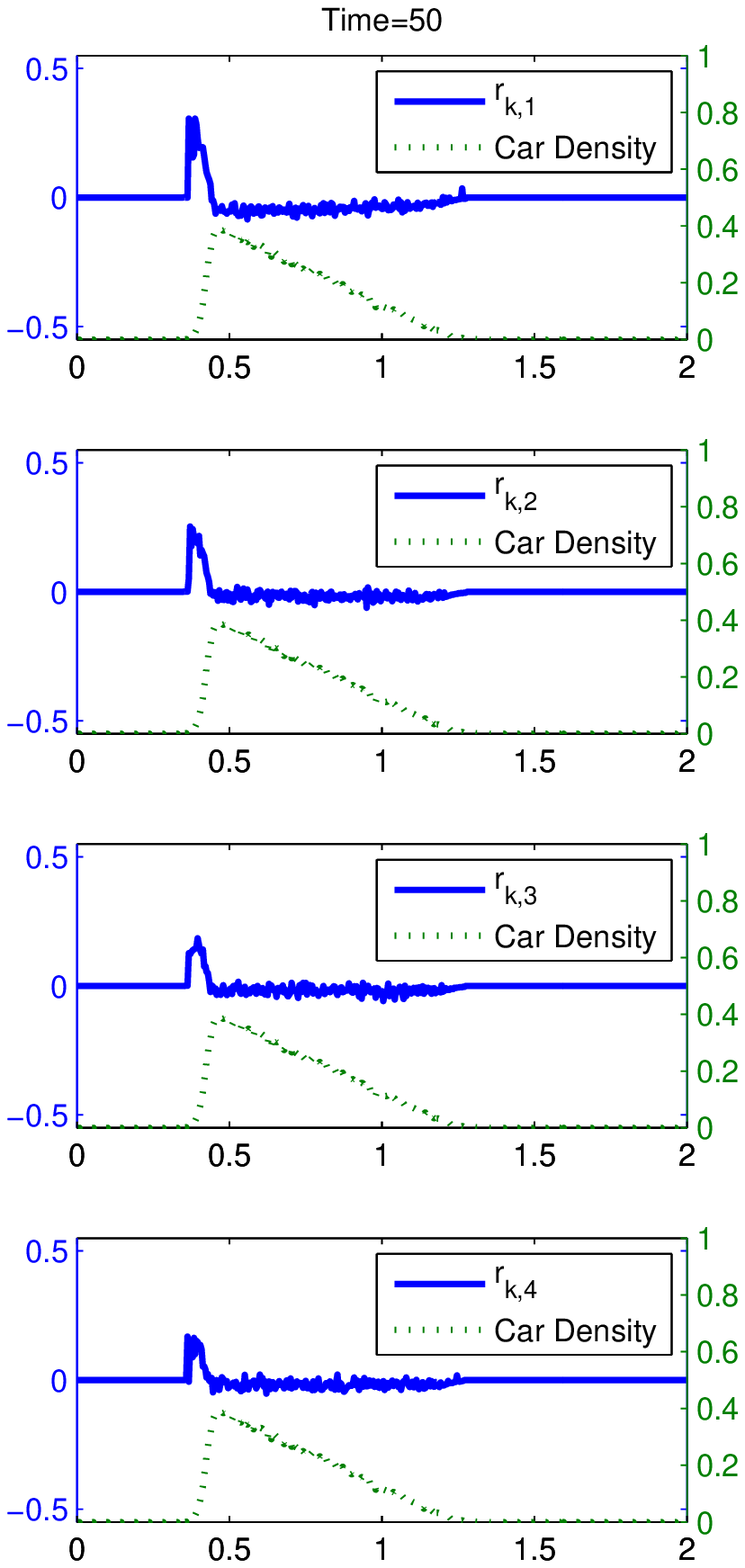}
\includegraphics[width=1.5in,height = 6.0in]{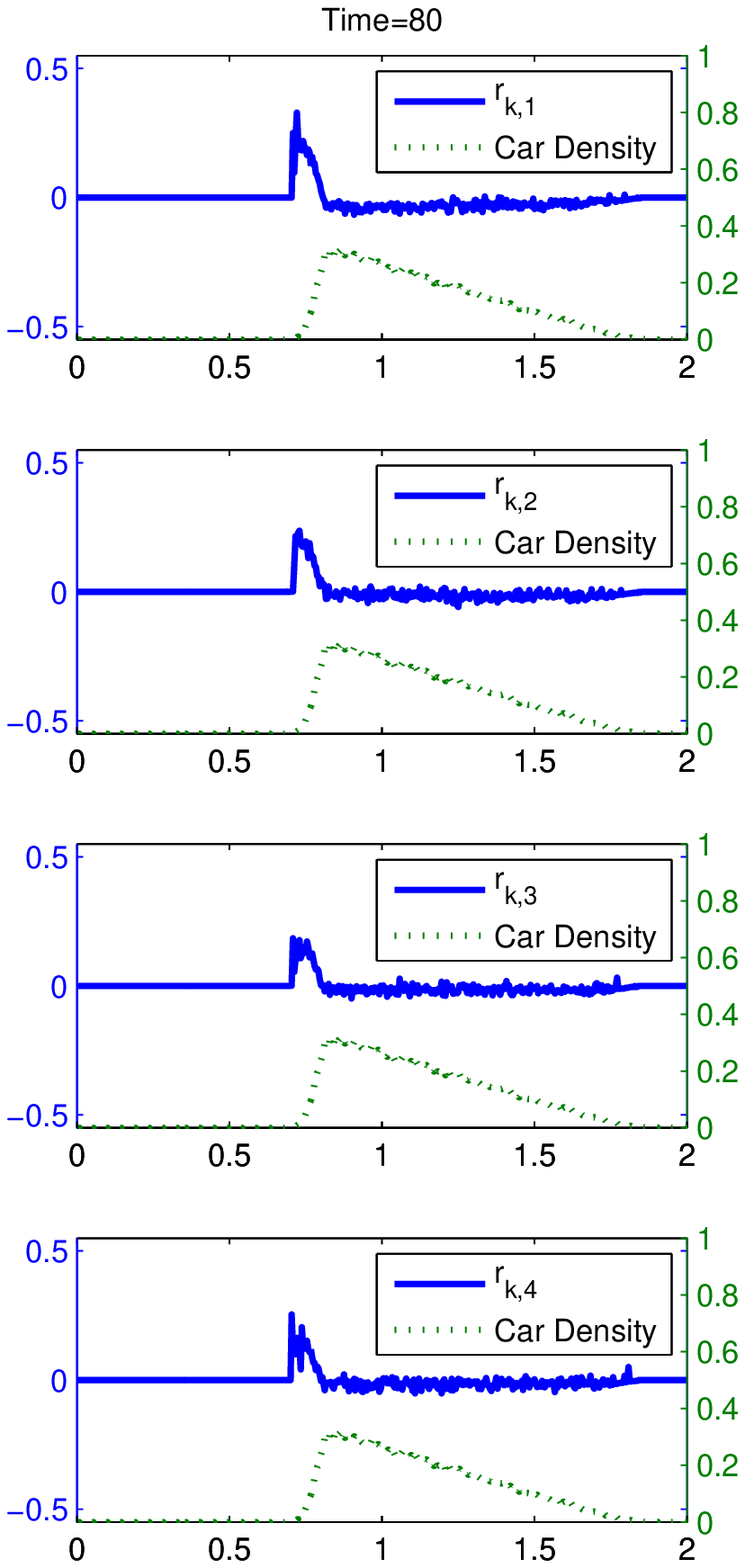}}
\caption{Left scale:  correlation coefficients $r_{k,k+p}$, $1 \le p \le 4$, for $\beta=0$
and $M = 0$, i.e., no look ahead potential. Right scale:  mean car density,}
\label{fig41}
\end{figure}
\end{center}

\begin{center}
\begin{figure}[H]
\centerline{
\includegraphics[width=1.5in,height = 6.0in]{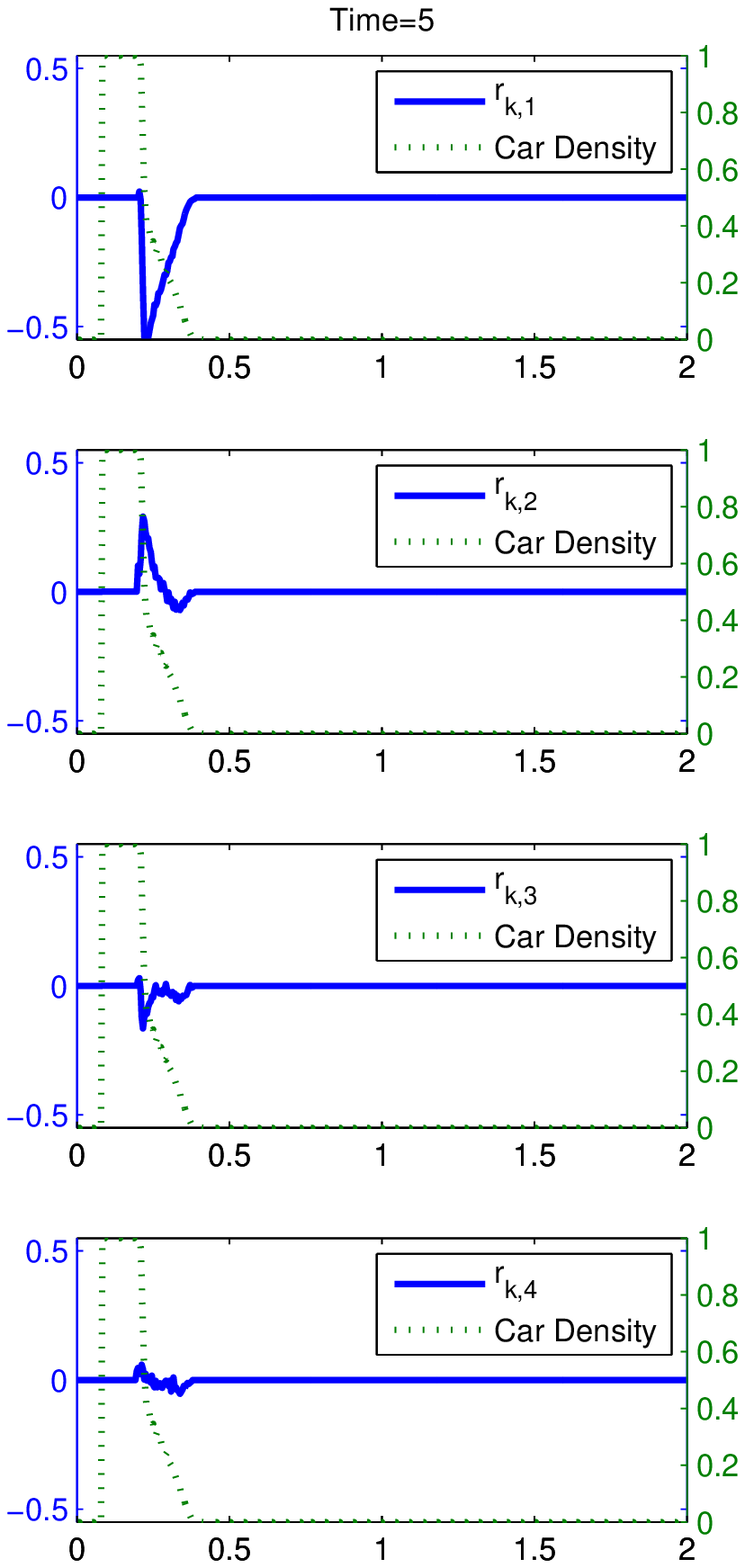}
\includegraphics[width=1.5in,height = 6.0in]{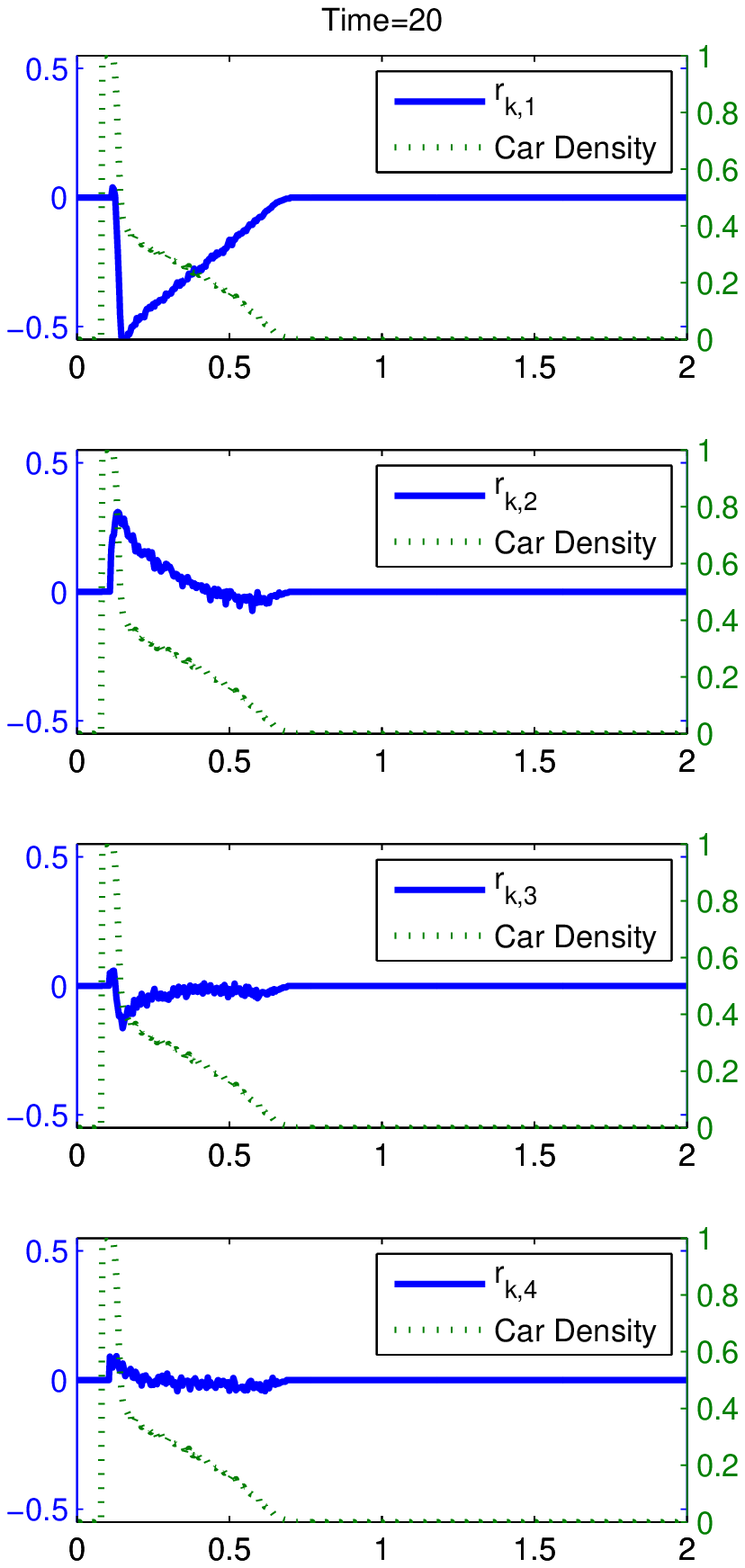}
\includegraphics[width=1.5in,height = 6.0in]{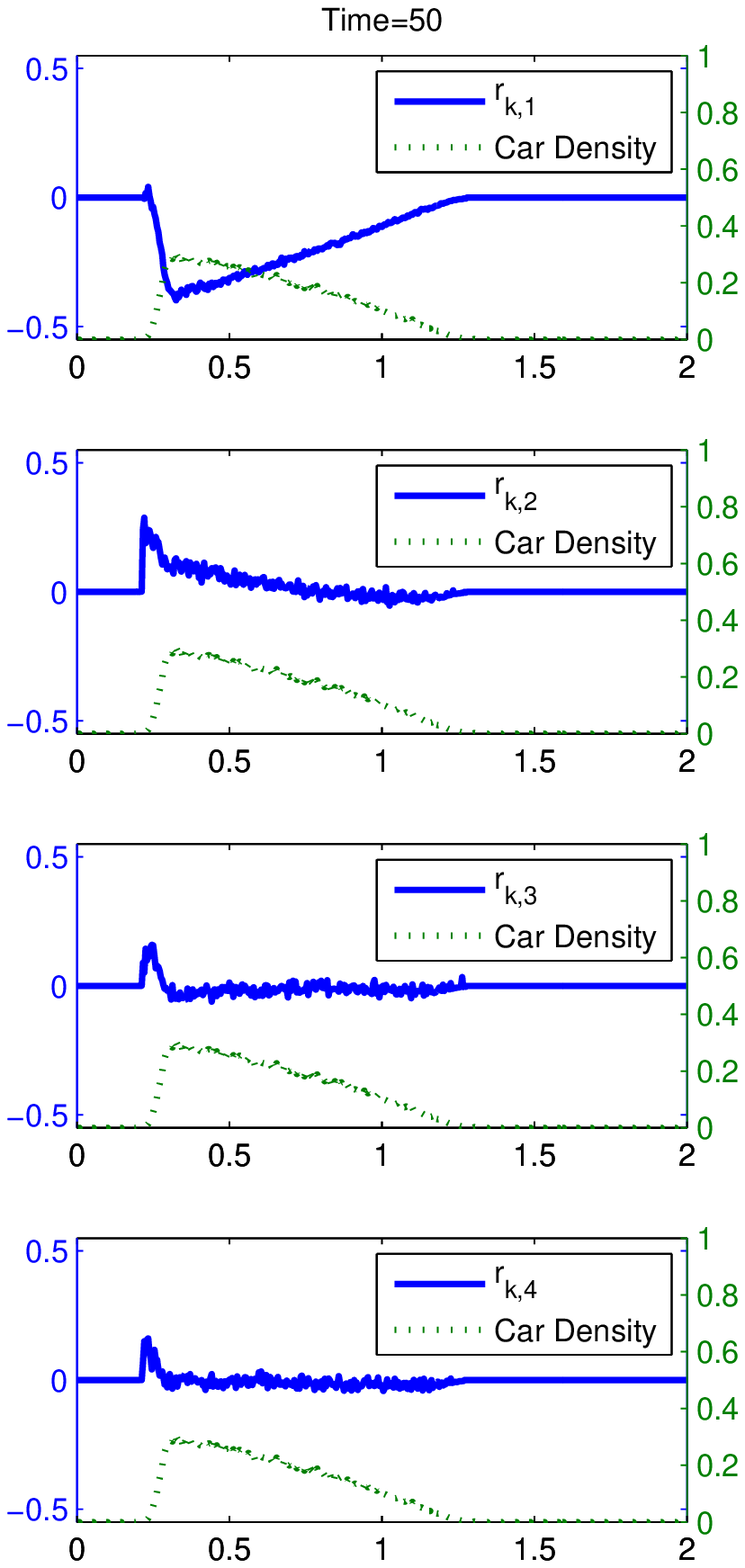}
\includegraphics[width=1.5in,height = 6.0in]{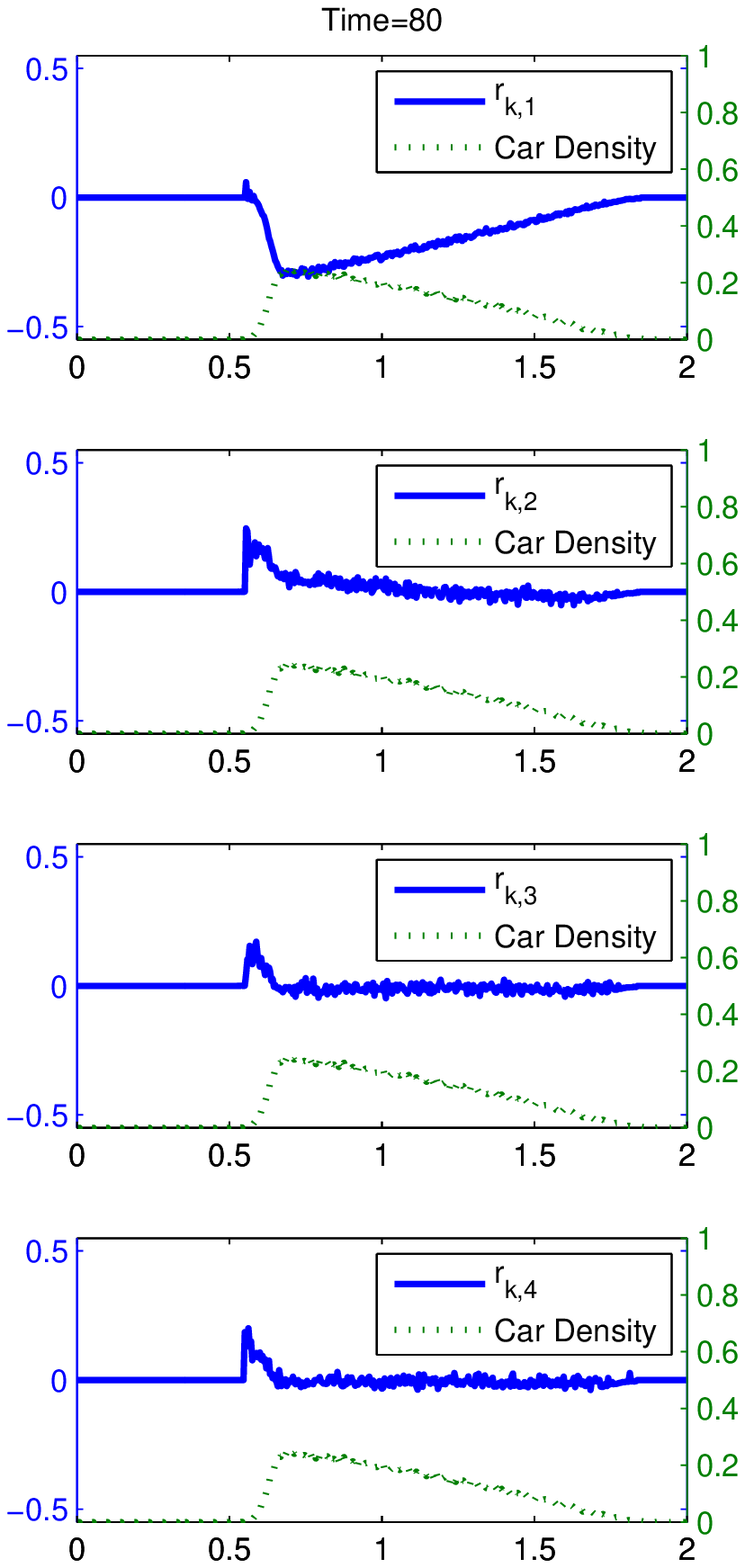}}
\caption{Left scale:  correlation coefficients  $r_{k,k+p}$, $1 \le p \le 4$, for $\beta=3$
and $M = 1$. Right scale:  mean car density,}
\label{fig42}
\end{figure}
\end{center}

\begin{center}
\begin{figure}[H]
\centerline{
\includegraphics[width=1.5in,height = 6.0in]{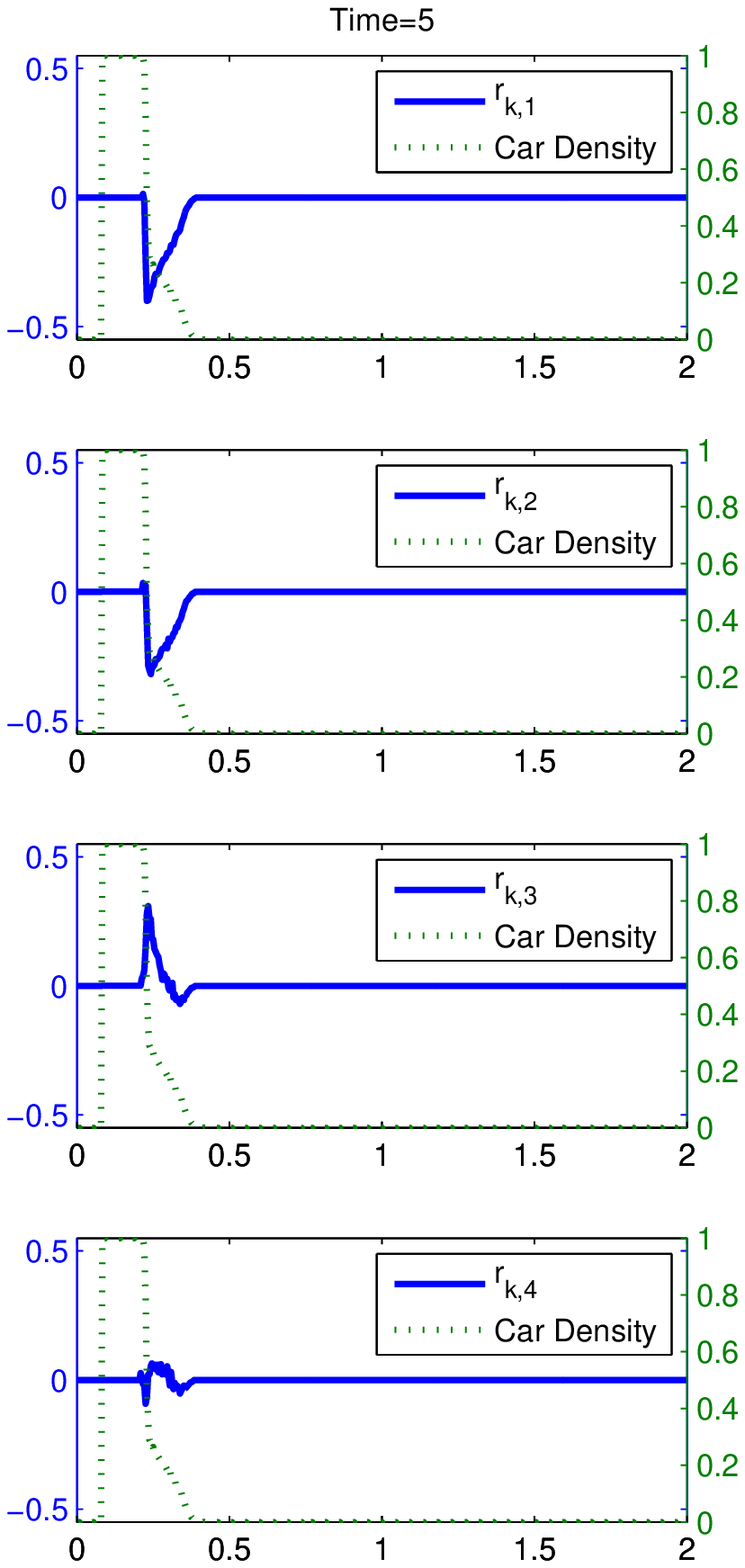}
\includegraphics[width=1.5in,height = 6.0in]{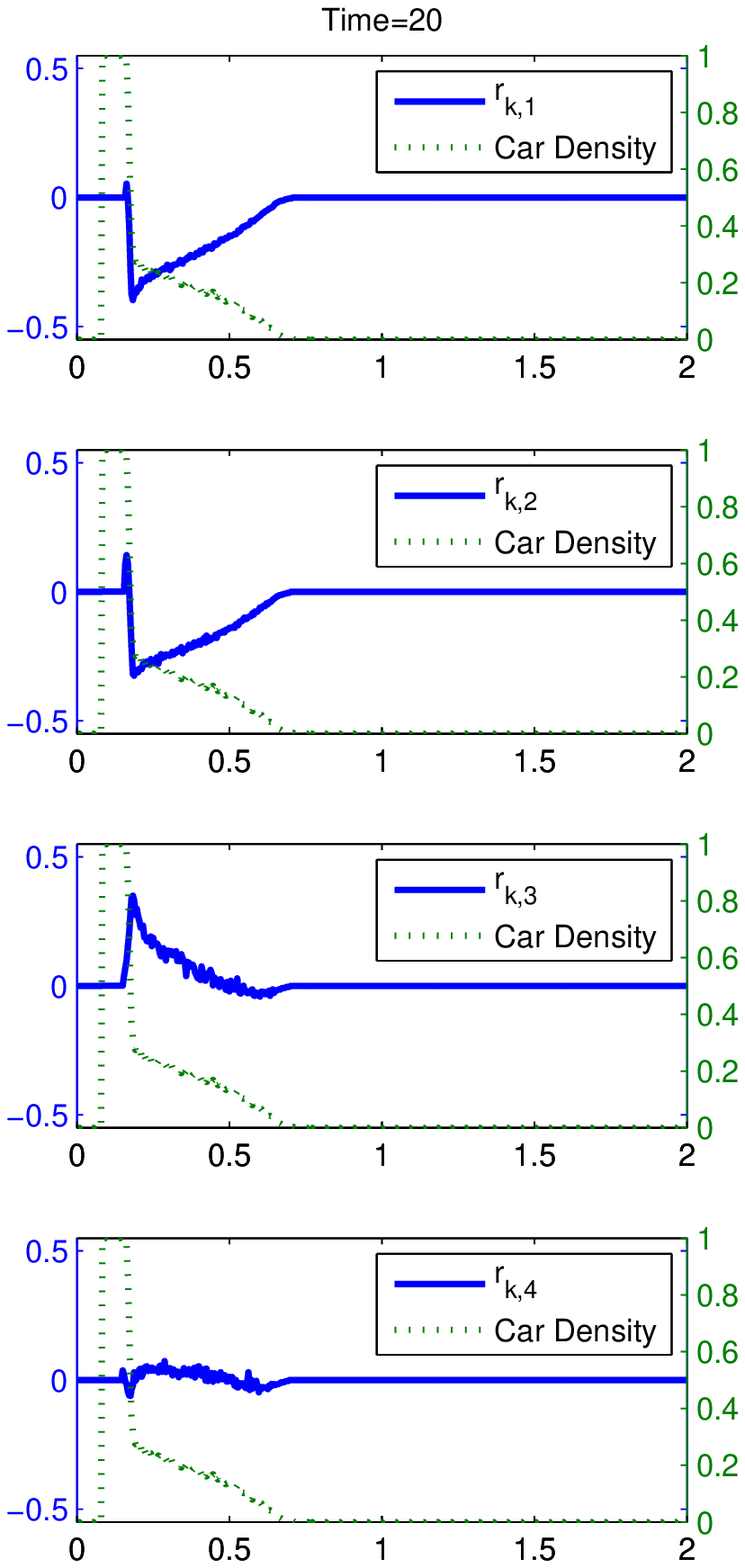}
\includegraphics[width=1.5in,height = 6.0in]{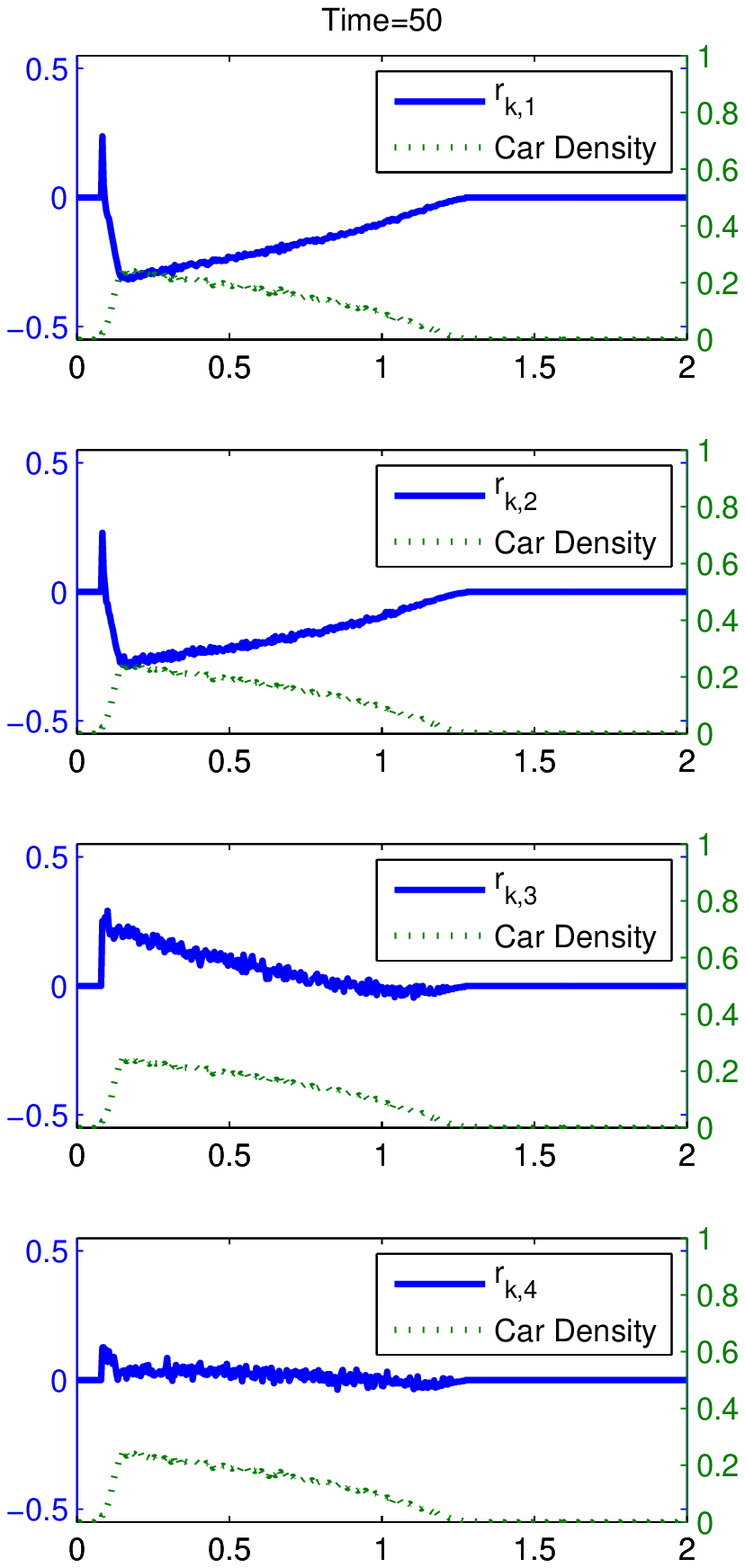}
\includegraphics[width=1.5in,height = 6.0in]{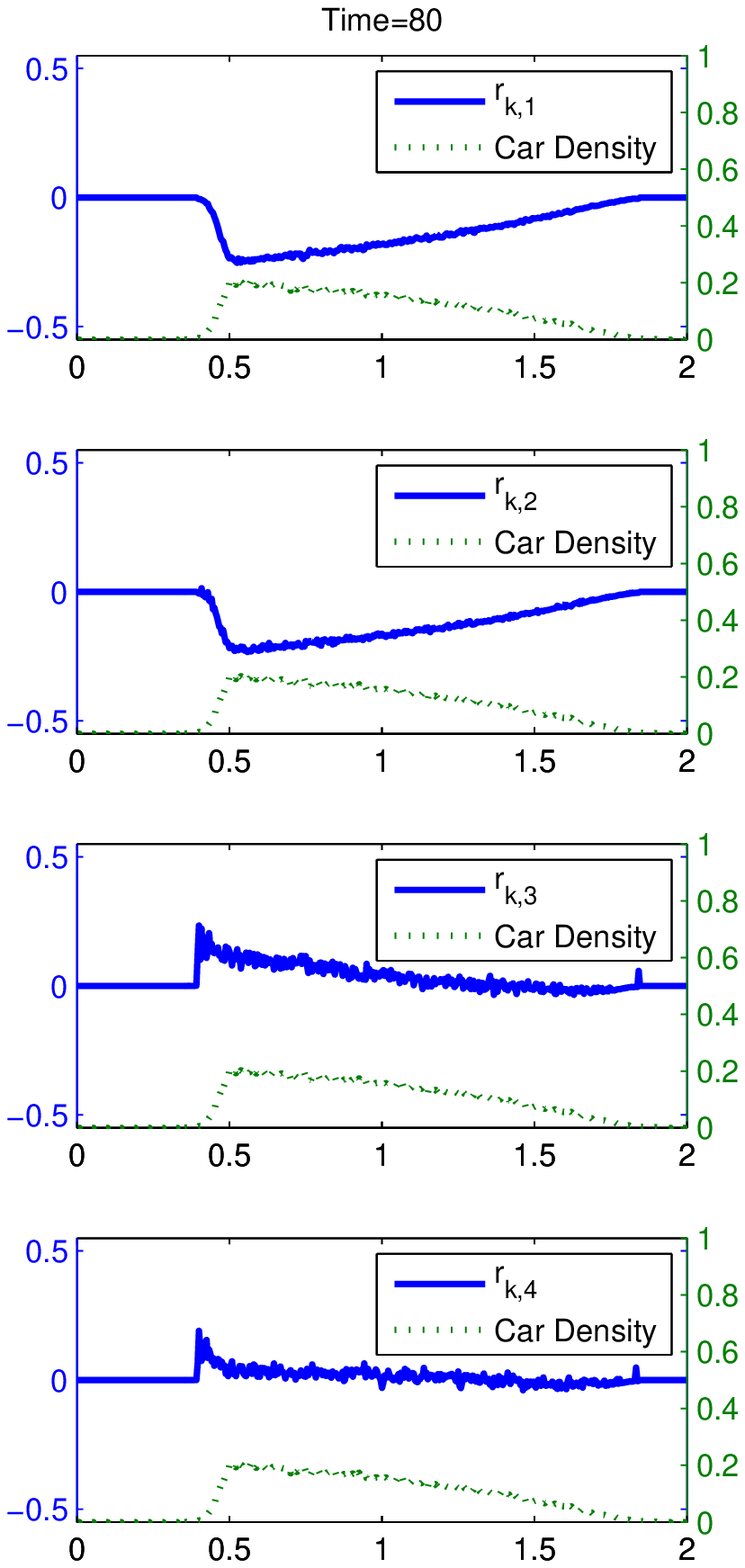}}
\caption{Left scale:  correlation coefficients  $r_{k,k+p}$, $1 \le p \le 4$, for $\beta=3$
and $M = 2$.  Right scale:  mean car density,}
\label{fig43}
\end{figure}
\end{center}

We compute correlation coefficients for $1
\le i \le 4$ using three sets of parameter values
\begin{enumerate}
     \item $\beta=0$ and $M = 0$, i.e, no look ahead potential (Figure
\ref{fig41});
     \item $\beta=3$ and $M = 1$ (Figure \ref{fig42});
     \item $\beta=3$ and $M = 2$ (Figure \ref{fig43}).
\end{enumerate}
The main conclusion are as follows:
\begin{itemize}
     \item 
     \red{
     Without the look ahead, the correlations are relatively weak, except
at the trailing front, which shows a moderately strong positive correlation (see
columns $2$--$4$ of Figure \ref{fig41}). 
At early time this correlation is stronger for smaller $i$, but comparable at
later times for $i >1$.  These correlations occur because any car that is slower than the ``main pack'' (maximum of the density in Monte-Carlo simulations)
will cause cars behind it to also slow down.  Since it is very likely that more than one car is slower than the main pack,
it produces positive correlations at the trailing front.
}
     \item \red{
     When $\beta$ is nonzero, the correlations $r_{k,k+i}$ depend strongly on the
value of $i$ relative to $M$.  For $i \leq M$, the correlations are large and negative, 
with the strongest correlation being at the trailing front.  This is due to the look-ahead potential 
which tends to maintain a spacing of $M$ zeros between cars.  For $i > M$, the correlations are similar in
the zero-potential case, especially at longer times.  
}
\end{itemize}

\subsection{Testing the Entire Closure for the New Model}

Next, we test the closure for the right-hand side of the mesoscopic model in
\eqref{den1} with the look-ahead dynamics, but just one car
in the look-ahead potential.  We select parameters $M=1$ and $\beta=3$, and
using sample statistics, test the closure
\begin{eqnarray}
&& \bE \left[
c_0 e^{-\beta \sigma_{k+1}} \sigma_{k-1} (1-\sigma_{k}) - 
c_0 e^{-\beta \sigma_{k+2}} \sigma_k (1-\sigma_{k+1})
\right] \approx
\nonumber \\
\label{rhsa1} \\
&& 
c_0 \left[1 + \bE \sigma_{k+1} (e^{-\beta} - 1)\right] \bE\sigma_{k-1}
(1-\bE\sigma_{k}) -
c_0 \left[1 + \bE \sigma_{k+2} (e^{-\beta} - 1)\right] \bE\sigma_k
(1-\bE\sigma_{k+1}).
\nonumber 
\end{eqnarray}
We compare these results to the closure which assumes no look-ahead 
potential, i.e., we test the closure
\begin{eqnarray}
&& \bE \left[
c_0 e^{-\beta \sigma_{k+1}} \sigma_{k-1} (1-\sigma_{k}) - 
c_0 e^{-\beta \sigma_{k+2}} \sigma_k (1-\sigma_{k+1})
\right] \approx
\nonumber \\
\label{rhsa2} \\
&&
c_0  \bE\sigma_{k-1} (1-\bE\sigma_{k}) - 
c_0  \bE\sigma_k (1-\bE\sigma_{k+1}).
\nonumber
\end{eqnarray}

\begin{center}
\begin{figure}[H]
\centerline{
\includegraphics[width=7.5in]{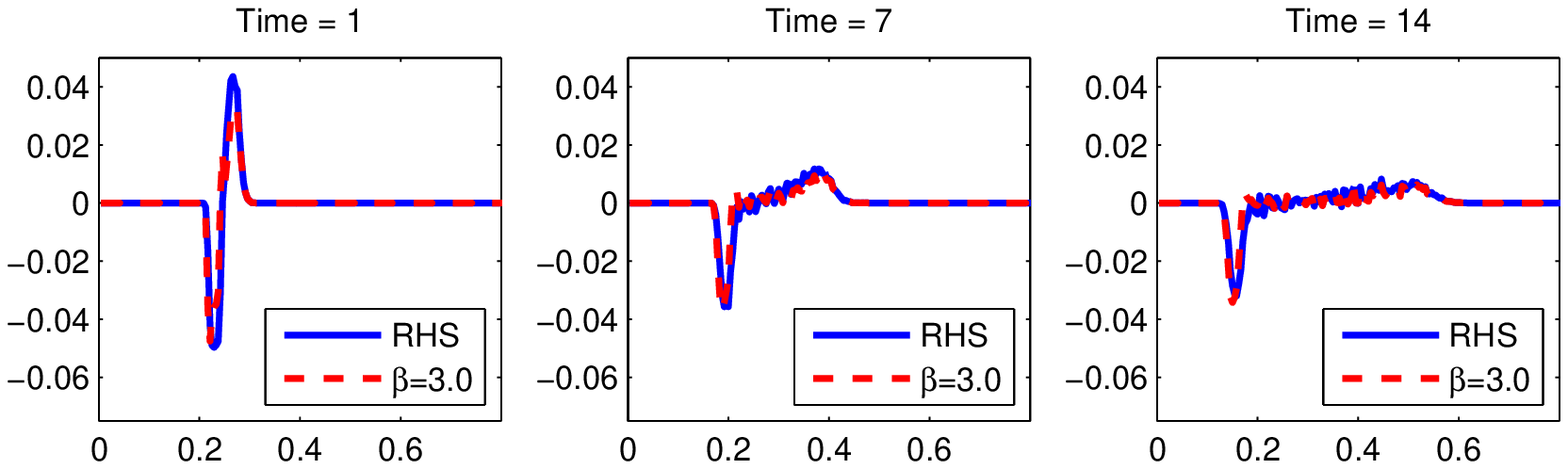}}
\centerline{
\includegraphics[width=7.5in]{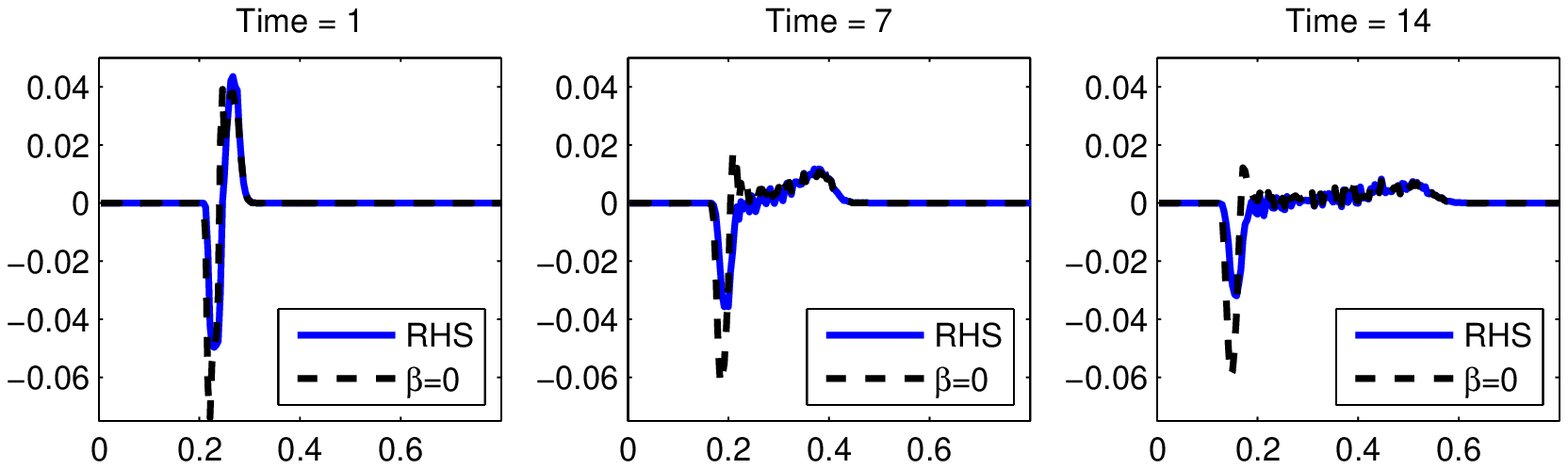}}
\centerline{
\includegraphics[width=7.5in]{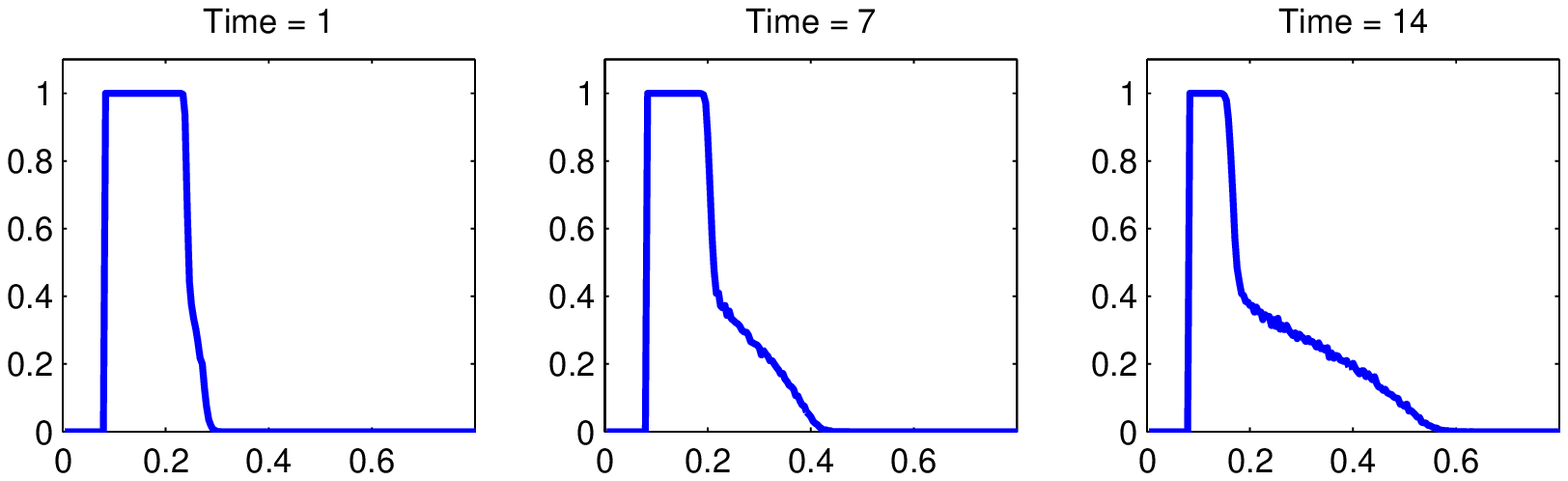}}
\caption{Comparison of closures from \eqref{rhsa1} (top line) and \eqref{rhsa2}
(middle line) during the course of a stochastic simulation.  Data is averaged
over five cell intervals to remove noise.  The bottom line is the density
profile, which is provided as a reference for locating the fronts.
The look ahead potential is $M=1$.}
\label{fig5}
\end{figure}
\end{center}

\begin{center}
\begin{figure}[H]
\centerline{
\includegraphics[width=7.5in]{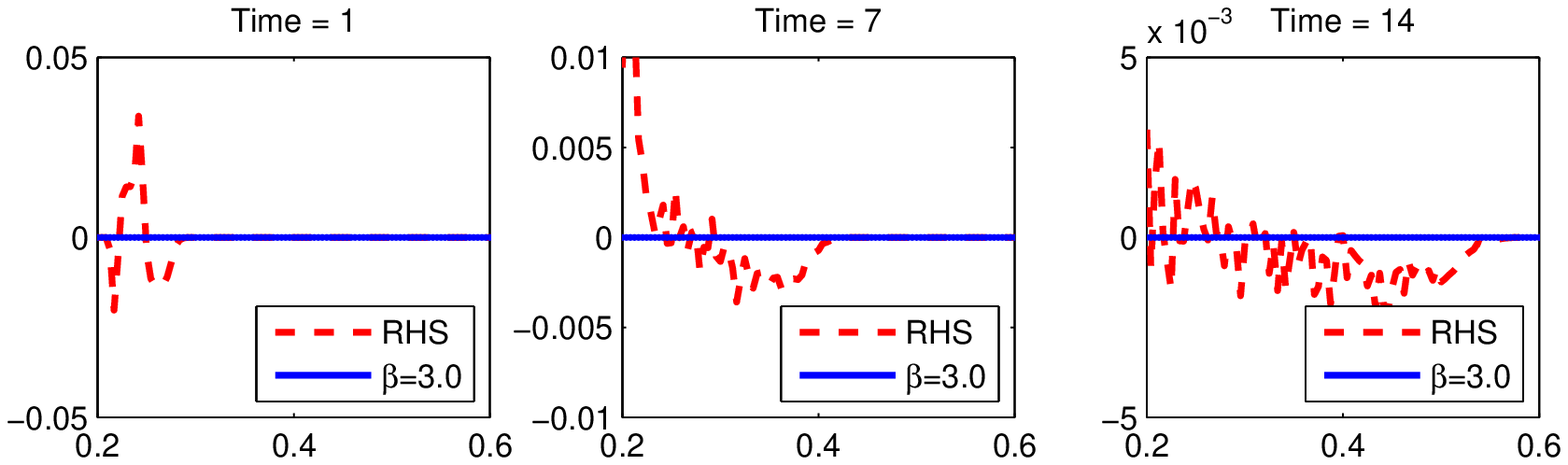}}
\centerline{
\includegraphics[width=7.5in]{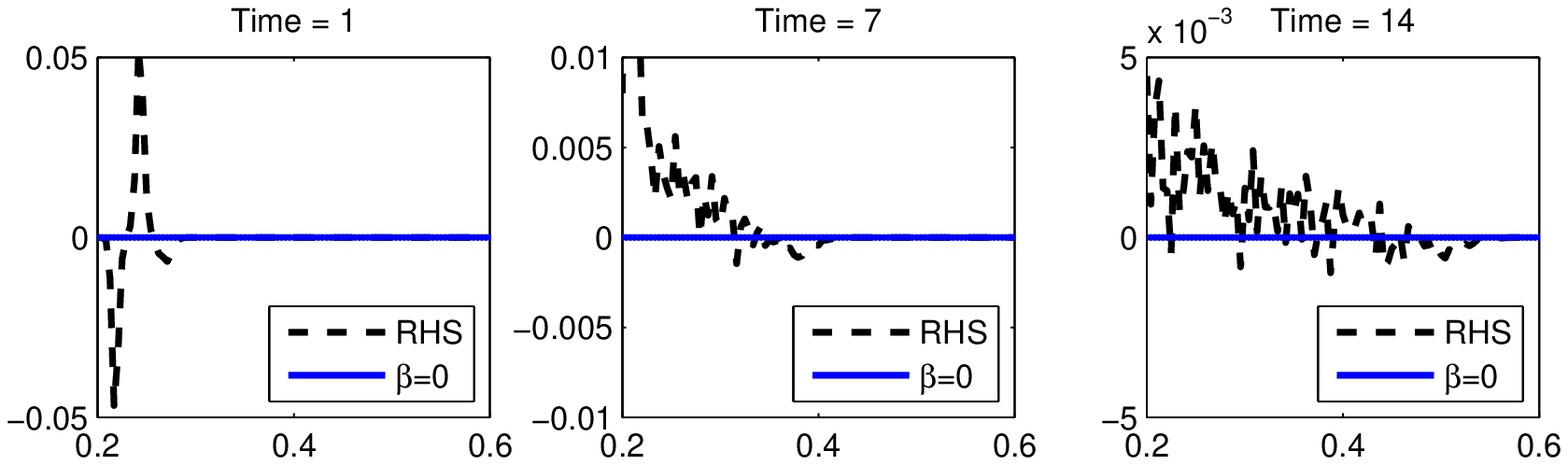}}
\caption{Comparison of closures from \eqref{rhsa1} (top line) and \eqref{rhsa2}
(middle line) during the course of a stochastic simulation.  Here we plot
the \textit{differences} between the exact closure and the two approximate
closures.  Thus positive values mean and overestimate of the time derivative.
Axes are magnified to allow inspection of the leading front.  Data is averaged
over five cell intervals to remove noise.
The look ahead potential is $M=1$.}
\label{fig6}
\end{figure}
\end{center}

\begin{center}
\begin{figure}[H]
\centerline{
\subfigure[$\beta=3$, $M=1$]{
\includegraphics[width = 7.0in]{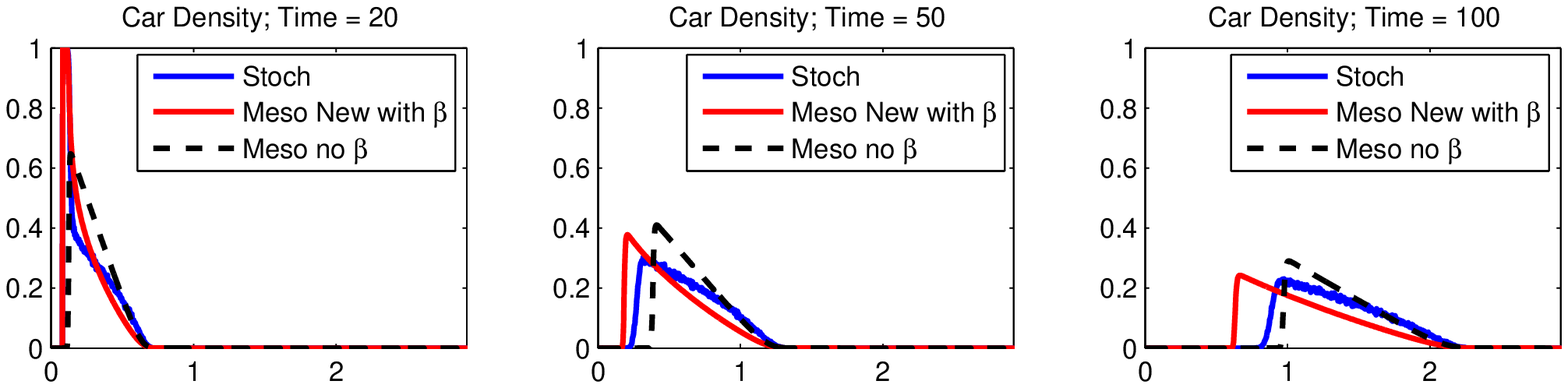}}}
\centerline{
\subfigure[$\beta=3$, $M=3$]{
\includegraphics[width = 7.0in]{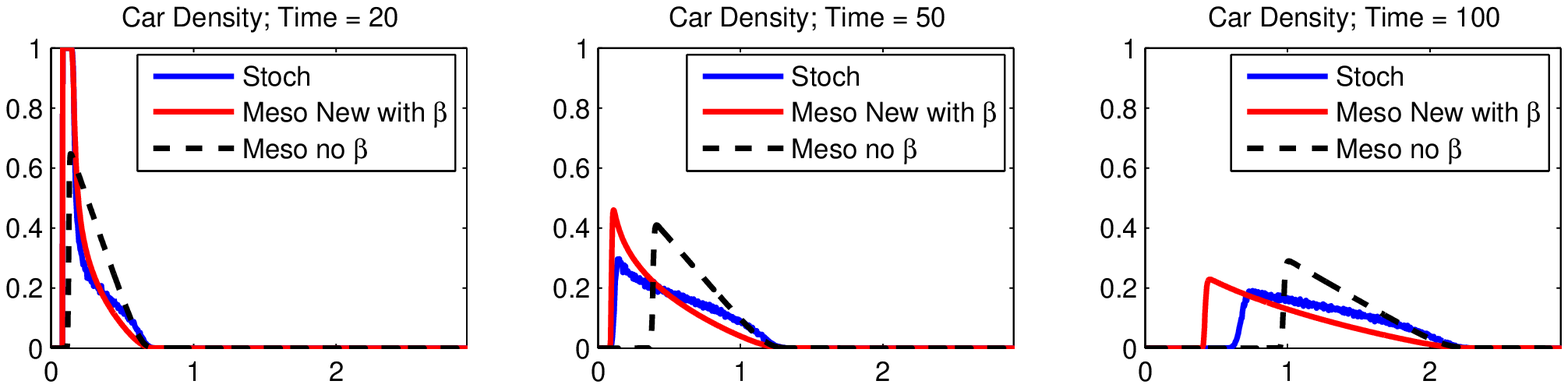}}}
\centerline{
\subfigure[$\beta=6$, $M=1$]{
\includegraphics[width = 7.0in]{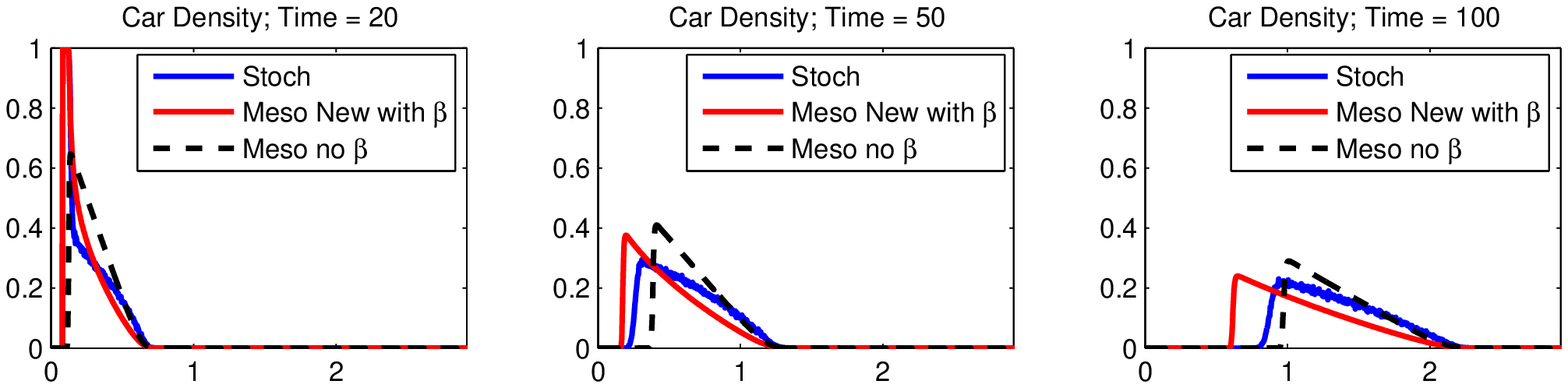}}}
\centerline{
\subfigure[$\beta=6$, $M=3$]{
\includegraphics[width = 7.0in]{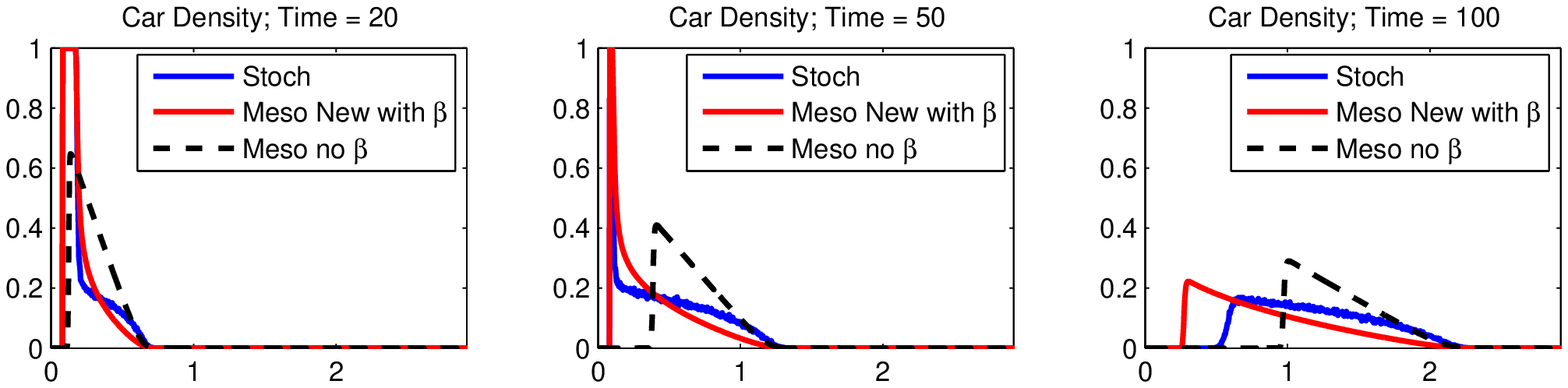}}}
\caption{Comparison of the car density in the Monte-Carlo simulations of the
stochastic model with $\beta=3$, $6$ $M=1$, $3$  and simulations of the new mesoscopic model
and the mesoscopic model without $\beta$ (i.e. $\beta=0$).}
\label{fig:pre_hack}
\end{figure}
\end{center}

Results from these experiments are given in Figures \ref{fig5} and \ref{fig6},
where we have averaged the expressions over five cells to reduce the noisy
fluctuations due to the stochastic nature of the equations. It
is clear in from Figure \ref{fig5} that the closure in \eqref{rhsa1} is 
more accurate near the trailing front.  However, a close inspection of the
leading front (see Figure \eqref{fig6}) shows that the closure in \eqref{rhsa2}
does a better job at the leading front.  Indeed the closure in \eqref{rhsa2}
underestimates the exact value; consequently not enough cars move
with the leading front.  This difference is quite subtle, but is observable
for times well beyond those give here.%
\footnote{The difference becomes lost in the noise around time $t=100$.}
Moreover, the effect in the density is quite noticeable.  From Figure
\eqref{fig:pre_hack}, one may observe the general trend that the closure in
\eqref{rhsa2} leads to better results at the leading front, while the
\eqref{rhsa1} performs better near the trailing front.

%
%
%
\section{Empirical Correction to the Equation for the Density}
\label{sec:hack}
As observed in the previous section, the proposed closure in \eqref{rhsa1}
underestimates the right-hand side of \eqref{den1}. This is also
evident from the plots in Figure \ref{fig:pre_hack} comparing the Monte-Carlo
simulations of the stochastic problem with simulations of the mesoscopic
equations. In particular, the bulk of cars behind the leading front not
propagate fast enough in the  mesoscopic model.  This means the effect of the
look-ahead dynamics in the proposed closure \eqref{rhsa1} is too strong.

On the other hand, the right-hand side in \eqref{rhsa2} (equivalent to setting
$\beta=0$ in the right-hand side) agrees much closure with the microscopic model. 
We therefore conjecture that the cross-correlations effectively
make the look-ahead parameter, $\beta$, in the mesoscopic model smaller.
For example in Figure \ref{fig:pre_hack}, the leading front in the solution of the
stochastic model is``in between'' the two solutions of the mesoscopic model with
and without $\beta$.

Based on these observations, we compensate for the error in the right-hand
side of \eqref{den1} by introducing a nonlinear 1ook-ahead parameter $\beta =
\beta \rho_k^d$, where $d$ is an empirical fitting parameter. Since $0 \le
\rho_k \le 1$ the proposed, nonlinear compensation of the look-ahead dynamics
leads to the effective reduction of the look-ahead potential. The resulting mesoscopic model is

\begin{align}
\frac{d}{dt} \rho_k &= 
c_0 \rho_{k-1} (1-\rho_k) 
\prod_{i=1}^{M} \left[ 1 + \rho_{k+i} \left( e^{-\beta' \rho_{k+i}^d} - 1
\right) \right]
\nonumber \\
\label{mesoemp} \\
& \quad -  c_0 \rho_k (1-\rho_{k+1})
\prod_{i=1}^{M} \left[ 1 + \rho_{k+i+1} \left( e^{-\beta' \rho_{k+i+1}^d} - 1
\right) \right].
\nonumber 
\end{align}

\begin{center}
\begin{figure}[H]
\centerline{
\subfigure[$\beta=3$, $M=1$, $d=2$]{
\includegraphics[width = 7.0in]{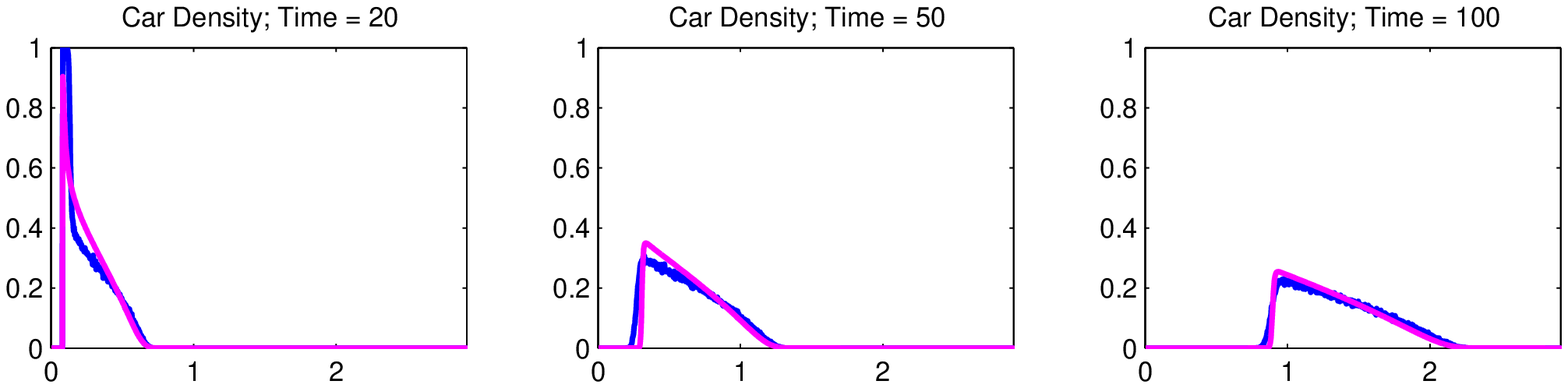}}}
\centerline{
\subfigure[$\beta=6$, $M=1$, $d=2$]{
\includegraphics[width = 7.0in]{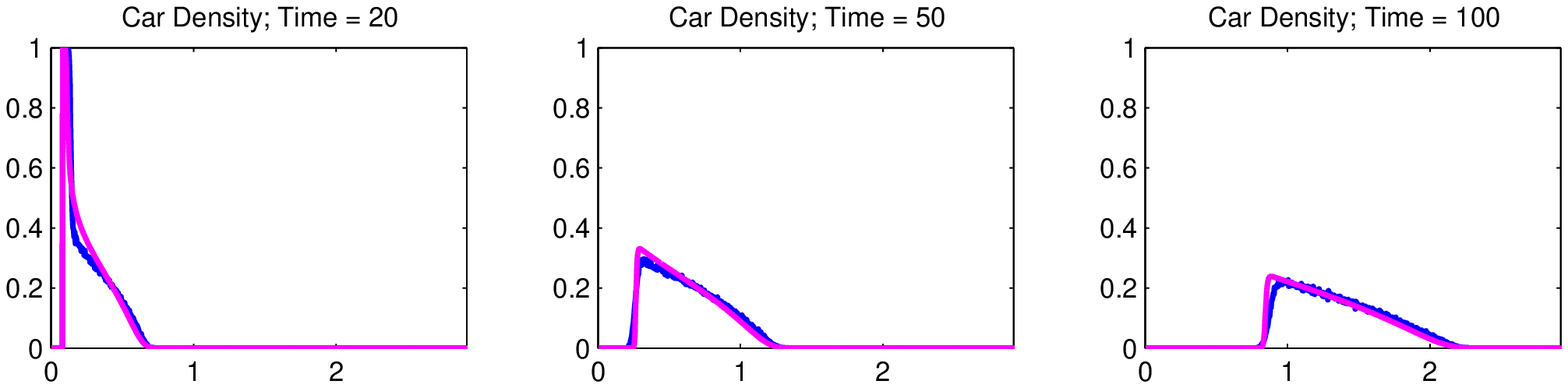}}}
\centerline{
\subfigure[$\beta=3$, $M=5$, $d=0.5$]{
\includegraphics[width = 7.0in]{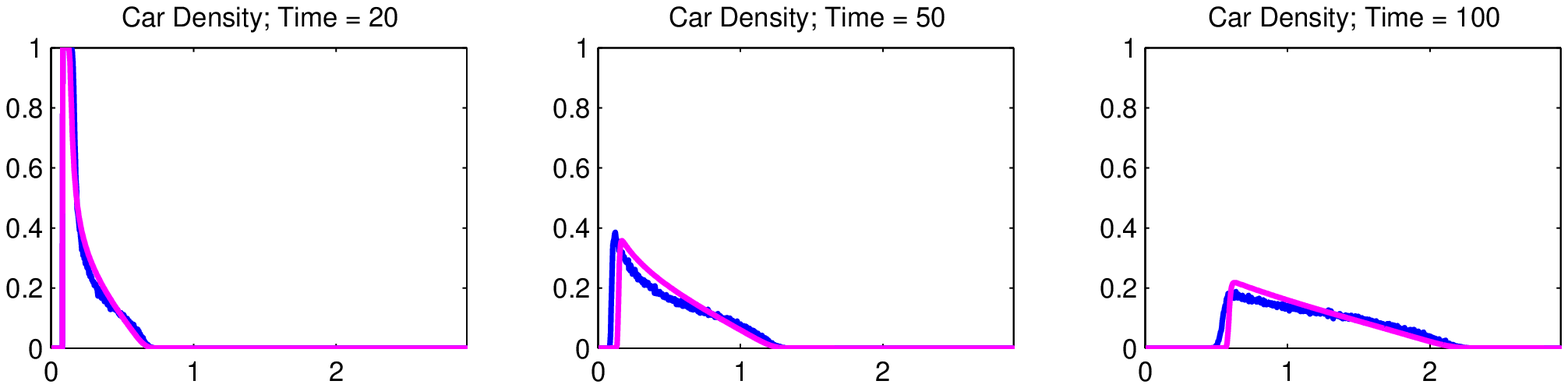}}}
\centerline{
\subfigure[$\beta=6$, $M=5$, $d=0.5$]{
\includegraphics[width = 7.0in]{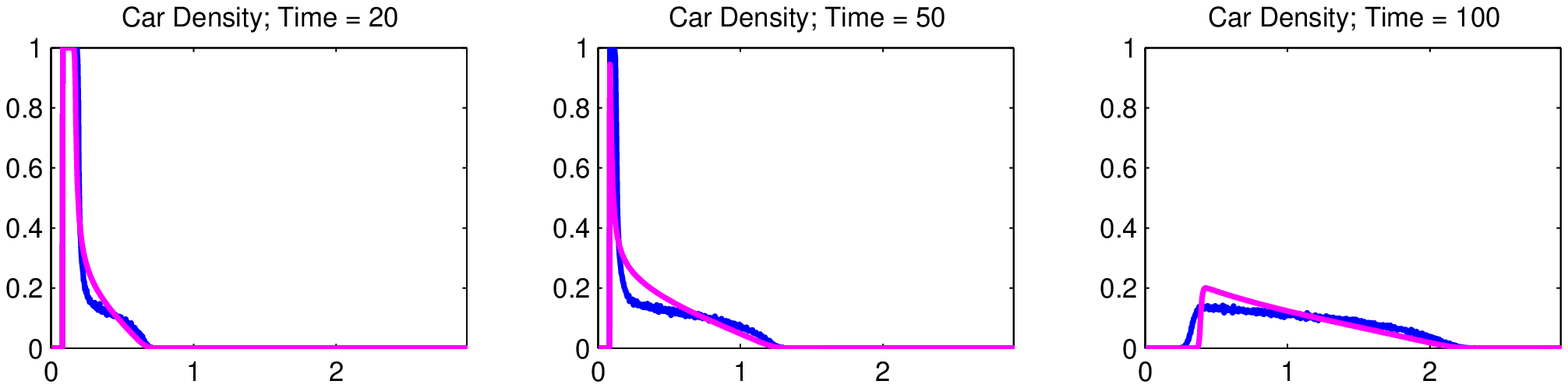}}}
\caption{Comparison of the car density in the Monte-Carlo simulations of the
stochastic model (blue line) and simulations of the mesoscopic model in
\eqref{mesoemp} (magenta). The exponent $d$ in \eqref{mesoemp} seems to depend 
strongly on $M$, but not on $\beta$.}
\label{fig:hack}
\end{figure}
\end{center}

Results comparing the Monte-Carlo simulations of the stochastic model and the
simulations of the modified mesoscopic model from \eqref{mesoemp} are depicted
in Figure \ref{fig:hack} for four different parameter combinations of $\beta$
and $M$.  As a reference, these results should be compared to those from Figure
\ref{fig:pre_hack}.  At this point, $d$ is a tuning parameter.  From experiments, we have found that a ``good'' choice of $d$
is sensitive to the value of $M$, but rather insensitive to the value $\beta$.

\subsection{Continuum Limits}

Following the arguments from Sections \ref{sec:meso_old} and
\ref{sec:meso_new}, it is straightforward to show that the continuum limit for
the empirical closure takes the form of the conservation law
\eqref{eq:cons_law}, but with a new flux given by
\begin{equation}
\phi(\bar{\rho})(x,t) = v_0 e^{-\frac{\beta}{L} \int_0^L \bar{\rho}^{d+1}(x+y,t)
dy}
\bar{\rho}(x,t)
(1-\bar{\rho}(x,t)) \:.
\end{equation}
As in \cite{Sopasakis-Katsoulakis-2006}, one may formally expand the integral
in the exponential.%
\footnote{
This is done for a general potential in
\cite{Sopasakis-Katsoulakis-2006};  here we consider only the constant
potential $\beta$.}  To second order, we find
\begin{equation}
     \frac{1}{L} \int_0^L \bar{\rho}^{d+1}(x+y,t) dy \approx \bar{\rho}^{d+1}(x)
+
\frac{L}{2} \frac{\partial}{\partial x} \bar{\rho}^{d+1}(x)
+ \frac{L^2}{6} \frac{\partial^2}{\partial x ^2}\bar{\rho}^{d+1}(x)
\end{equation}
so that, via Taylor expansion,
\begin{equation}
     e^{-\frac{\beta}{L} \int_0^L \bar{\rho}^{d+1}(x+y,t) dy} \approx
     e^{-\beta \bar{\rho}^{d+1}(x)}
     \left(1 - \frac{\beta L}{2} \frac{\partial}{\partial x} \bar{\rho}^{d+1}(x)
 - \frac{\beta L^2}{6} \frac{\partial^2}{\partial x ^2}\bar{\rho}^{d+1}(x) \right)
\end{equation}
This gives the following local, continuum model
\begin{equation}
\bar{\rho}_t
+  \left( e^{-\beta  \bar{\rho}^{d+1}} \bar{\rho} (1 - \bar{\rho}) \right)_x
= 
 \frac{\beta L}{2} \left(  \bar{\rho} (1 - \bar{\rho}) e^{-\beta \bar{\rho}^{d+1}}
 (\bar{\rho}^{d+1})_x \right)_x
+
\frac{\beta L^2}{6} \left(  \bar{\rho} (1 - \bar{\rho}) e^{-\beta  
\bar{\rho}^{d+1}}
(\bar{\rho}^{d+1})_{xx} \right)_x
\end{equation}
which introduces an additional layer of nonlinearity to the models presented in
\cite{Sopasakis-Katsoulakis-2006}.  We note, however, that the proceeding
derivation is completely formal.  In particular, it assumes that solutions
are smooth and that the error in the Taylor series expansions are small.

%
%
\section{Discussion and Conclusions}
\label{sec:conclusions}

In the paper, we have examined in detail the cellular automata (CA) traffic
model introduced in \cite{Sopasakis-Katsoulakis-2006}, including a host of
numerical experiments performed on a rarefaction wave evolving from a
``red-light'' initial condition.  The work in \cite{Sopasakis-Katsoulakis-2006}
included mesoscopic (ODE) and continuum (PDE) approximations of the cellular
automata model.  For the mesoscopic model, two assumptions are required.  The
first of these has been removed in the current paper, leading to an improved
model at the ODE level, but with the same PDE in the continuum limit. 
Since the look-ahead potential conidered here have also been used in other contexts
\cite{kkpv2012,paka2012,kmk2003}, the improved mesoscopic model can lead to 
significant improvements in other areas.

The second assumption used in \cite{Sopasakis-Katsoulakis-2006} to derive the
mesoscopic model is based on the independence of spatial neighbors in
the stochastic process defined by the CA model. We have shown numerically, that
such an assumption does not hold for strong look-ahead potentials and moreover,
that the correlations effectively weaken the effect of the look ahead
potential at the macroscopic level.  We conjecture that this effect is density dependent
and introduce an ad-hoc macroscopic
parameter which scales like a power law in the density.  The
choice of the exponent in the power law was determined experimentally and found
to depend strongly on the look ahead distance $M$.  However, for a given $M$,
the exponent behaves quite well for different $\beta$ and over long time
scales.  Finally, using this new ad-hoc model, we have derived nonlinear
variants of the known, local continuum models derived from the non-local
continuum model in \cite{Sopasakis-Katsoulakis-2006}. (See also Section
\ref{sec:meso_old}.)

We see three avenues for future work.  
First is the need to explore the
generality of the numerical results presented here.  Certainly a complete
repeat of the experiments in the current paper should be done for a
jam---that is, the movement of faster cars at low density into a
slower, higher-density regime.  For a continuum hyperbolic model, this would
correspond to a shock.  
The second issue to explore is the relationship between
the exponent $d$ and the look ahead distance $M$.  While the existence of such
a relationship can be inferred by our numerical experiments, a formula would
prove quite useful in practice. Finally, for engineering
purposes like traffic microscopic models are not computationally
practical. A possible alternative is to use filtering
\cite{Anderson-Moore-2011}, where the macroscopic or mesoscopic models are
corrected in real-time using measurements.  Such an approach has already been
investigated, for example, in \cite{mibo04,wapa05,smh03}.

\bigskip

{\bf Acknowledgements.}
The work presented in this paper emerged as a result of discussions in  
a working group at the NSF funded Statistical and Applied Mathematical  
Sciences Institute.  I. Timofeyev also acknowledges support from SAMSI as a  
long-term visitor in the Fall of 2010 and the Fall of 2011.

\bibliographystyle{siam}
\bibliography{traffic}
\end{document}